# Embedding with a Rigid Substructure

Igor Najfeld

Timothy F. Havel

*Biological Chemistry and Molecular Pharmacology*

*Harvard Medical School, Boston, MA 02115*

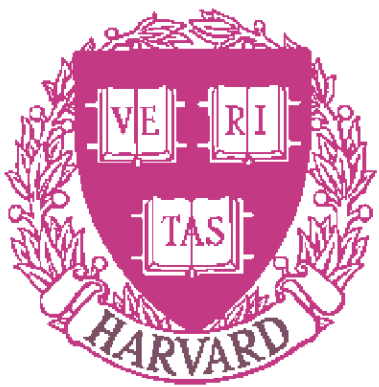

# ABSTRACT

This paper presents a new distance geometry algorithm for calculating atomic coordinates from estimates of the interatomic distances, which maintains the positions of the atoms in a known rigid substructure. Given an $M \times 3$ matrix of coordinates for the rigid substructure $\mathbf{X}$, this problem consists of finding the $N \times 3$ matrix $\mathbf{Y}$ that yields of global minimum of the so-called *STRAIN*, i.e.

$$\min_{\mathbf{Y}} \left\| \begin{bmatrix} \mathbf{X}\mathbf{X}^\mathsf{T} & \mathbf{X}\mathbf{Y}^\mathsf{T} \\ \mathbf{Y}\mathbf{X}^\mathsf{T} & \mathbf{Y}\mathbf{Y}^\mathsf{T} \end{bmatrix} - \begin{bmatrix} \mathbf{A} & \mathbf{B} \\ \mathbf{B}^\mathsf{T} & \mathbf{C} \end{bmatrix} \right\|_\mathsf{F}^2 ,$$

where $\mathbf{A} = \mathbf{X}\mathbf{X}^\mathsf{T}$, and $\mathbf{B}$, $\mathbf{C}$ are matrices of inner products calculated from the estimated distances.

The vanishing of the gradient of the *STRAIN* is shown to be equivalent to a system of only six nonlinear equations in six unknowns for the inertial tensor associated with the solution $\mathbf{Y}$. The entire solution space is characterized in terms of the geometry of the intersection curves between the unit sphere and certain variable ellipsoids. Upon deriving tight bilateral bounds on the moments of inertia of any possible solution, we construct a search procedure that reliably locates the global minimum. The effectiveness of this method is demonstrated on realistic simulated and chemical test problems.

# 1. Introduction

Distance geometry is widely used in the generation of three-dimensional atomic coordinates from chemical diagrams [1], in the determination of conformation from NMR data [2], in the modeling of protein structure by homology [3], and in docking ligands in their receptor's binding site [4]. More generally, distance geometry is a formalism within which a large variety of conformational problems can be stated in purely geometric terms [5]. Solving these problems generally involves finding three-dimensional coordinates that satisfy given lower and upper bounds on the interatomic distances. In contrast to relatively complex models of the intramolecular potential energy, this simple geometric problem formulation makes it possible to use rigorous mathematical techniques to assist the search for globally optimal solutions.

The most important of these mathematical techniques is the *EMBED* algorithm, which rapidly finds coordinates that closely fit the distance bounds [6]. As explained in the next section, these coordinates are the global minimum of a certain matrix optimization problem, which can be reliably located by eigenvalue methods. The availability of such good starting coordinates greatly improves the efficiency and con-vergence properties of their subsequent nonlinear optimization versus an *error function*, which measures the violations of the distance bounds [7]. One shortcoming of the *EMBED* algorithm is that there is no good way to weight the distances according to their precision, or exactly realize rigid substructures even when they are known in advance. This shortcoming is particularly serious in problems with large rigid substructures, as in the docking of a ligand to a receptor protein of known structure, or in building the loops and sidechains of a protein onto a backbone obtained via homology modeling.

In this paper we formulate a new minimization criterion which fixes the positions of the atoms in a known rigid substructure. The overall method may be regarded as a block matrix extension of the *EMBED* algorithm. We derive the necessary conditions for an optimal solution in the form of a nonlinear matrix equation. By considering the inertial tensor associated with a solution matrix, we have transformed this matrix equation into a system of six nonlinear equations among three eigenvalues and three eigenvectors of the inertial tensor. The solutions to this reduced nonlinear system are characterized in terms of the intersection curves between variable ellipsoids and the unit sphere: each eigenvector lies on the intersection curve associated with a given eigenvalue.

Following the establishment of the intersection conditions, we derive sharp bounds on a six-dimensional subspace containing all possible inertial tensors of interest, and in particular the inertial tensor of the global minimum. These bounds, together with an explicit parametrization of the intersection curves and the symmetries in the equations, are used to construct a specialized search procedure that systematically samples this subspace in order to extract a set of approximate solutions to the equations. We then refine this set of approximate solutions by means of iterative methods specifically designed for this task. The outcome is, in general, a rather small set of distinct solutions from which the global minimum is found by simply comparing their function values. The entire algorithm is evaluated on a realistic set of simulated and chemical test problems, which demonstrate its effectiveness in solving the above-mentioned conformational problems.

## 2. Background and Problem Formulation

We begin with some necessary background on the *EMBED* algorithm. Suppose we have estimates of all the interatomic distances $[d_{ij}]$

in a structure, and we wish to find a matrix of three-dimensional coordinates $\mathbf{X} = [\mathbf{x}_1, \ldots, \mathbf{x}_M]^\mathsf{T} \in \mathbf{R}^{M \times 3}$ for the $M$ atoms such that the distances calculated from the coordinates $\|\mathbf{x}_i - \mathbf{x}_j\|$ in some sense closely match their estimated values. The obvious way to do this is to minimize the weighted *STRESS* function

$$\sum_{1 = i < j}^{M, M} (w_{ij}[\|\mathbf{x}_i - \mathbf{x}_j\| - d_{ij}])^2 , \tag{1}$$

or *SSTRESS* function

$$\sum_{1 = i < j}^{M, M} (w_{ij}[\|\mathbf{x}_i - \mathbf{x}_j\|^2 - D_{ij}])^2 , \tag{2}$$

with respect to the coordinates, where $D_{ij} \equiv d_{ij}^2$ and $w_{ij} \geq 0$ are weights.

Although this can certainly be done [8], it turns out that another function, called the *STRAIN*, has nicer mathematical properties. To explain this, we expand eq. (2) as

$$4 \sum_{1 = i < j}^{M, M} w_{ij}^2 \, [(\mathbf{x}_i \cdot \mathbf{x}_j) - (\|\mathbf{x}_i\|^2 + \|\mathbf{x}_j\|^2 - D_{ij})/2]^2 . \tag{3}$$

By the law of cosines, this can be regarded as four times the weighted sum of the squares of the differences between the dot product $\mathbf{x}_i \cdot \mathbf{x}_j$ and an estimate thereof, namely $(\|\mathbf{x}_i\|^2 + \|\mathbf{x}_j\|^2 - D_{ij})/2$. This estimate is inconvenient since it depends on the coordinates we are trying to calculate. An estimate that is independent of the coordinates may be obtained from the following formula, which gives the squared distances to the center of mass of the configuration, $\overline{M} \equiv \sum_j m_j$, in terms of the squared distances among the atoms:

$$D_{0i} = \overline{M}^{-1} \sum_{j=1}^{M} m_j D_{ij} - \overline{M}^{-2} \sum_{1 = j < k}^{M, M} m_j m_k D_{jk} \ . \tag{4}$$

It is straightforward to show that when $D_{ij} = \|\mathbf{x}_i - \mathbf{x}_j\|^2$ for some set of center of mass coordinates, the estimate $D_{0i}$ is exact, i.e. $D_{0i} = \|\mathbf{x}_i\|^2$ [6]. Moreover, because eq. (4) is an averaging procedure, the $D_{0i}$ obtained from it are quite insensitive to errors in the $D_{ij}$.

The complete matrix of dot products between all pairs of coordinate vectors, or *Gram matrix*, is readily calculated from the coordinates as $\mathbf{X}\mathbf{X}^\top = [\mathbf{x}_i \cdot \mathbf{x}_j]_{i,j=1}^{M,M}$. If the coordinates $\mathbf{X}$ are center of mass coordinates and $\mathbf{D} = [D_{ij}]_{i,j=1}^{M,M} = [\|\mathbf{x}_i - \mathbf{x}_j\|^2]_{i,j=1}^{M,M}$, then it can also be written as $\mathbf{X}\mathbf{X}^\top = [(D_{0i} + D_{0j} - D_{ij})/2]_{i,j=1}^{M,M}$. Let $\mathbf{I}$ be an identity matrix, $\mathbf{1}$ a column vector of ones, and $\mathbf{m} = [m_1, \ldots, m_M]^\top$. Then this transforma-tion from matrices of squared distances to Gram matrices can be written succinctly (see [9]) as

$$\mathbf{X}\mathbf{X}^\mathsf{T} = -\frac{1}{2} [\mathbf{I} - (\mathbf{1}\mathbf{m}^\mathsf{T})/\overline{M}] \, \mathbf{D} \, [\mathbf{I} - (\mathbf{m}\mathbf{1}^\mathsf{T})/\overline{M}] \, , \qquad (5)$$

The mathematical equivalent of the Gram matrix in statistics is the well known covariance matrix.

If the $D_{ij}$ are only inexact estimates of the actual squared distances $\|\mathbf{x}_i - \mathbf{x}_j\|^2$, we can calculate an estimate of the corresponding Gram matrix as $\mathbf{A} = [a_{ij}]_{i,j=1}^{M,M} = [(D_{0i} + D_{0j} - D_{ij})/2]_{i,j=1}^{M,M}$. The *STRAIN* may then be defined as

$$\sum_{i,j=1}^{M,M} (w_{ij}[(\mathbf{x}_i \cdot \mathbf{x}_j) - a_{ij}])^2 = \|\mathbf{W} \bullet [\mathbf{X}\mathbf{X}^\mathsf{T} - \mathbf{A}]\|_\mathsf{F}^2 \, , \qquad (6)$$

where $\|\mathbf{C}\|_\mathsf{F}^2 = \sum_{ij} c_{ij}^2$ denotes the squared *Frobenius norm* of matrix $\mathbf{C}$, while $\mathbf{W} = [w_{ij}]_{i,j=1}^{M,M}$ is a matrix of weights and the "•" denotes the *Hadamard* (or entry-by-entry) matrix product.

This problem of minimizing the *STRAIN* does not have a simple solution for arbitrary weights, but if $w_{ij} = m_i m_j$ for some assignment of masses to the atoms, eq. (6) can be written as $\|\mathbf{M}[\mathbf{XX}^\mathsf{T} - \mathbf{A}]\mathbf{M}\|_\mathsf{F}^2$, where $\mathbf{M} = \mathit{diag}(\mathbf{m})$ is a nonsingular diagonal matrix containing the masses. The mass-weighted Gram matrix $\mathbf{M}[\mathbf{XX}^\mathsf{T}]\mathbf{M}$ is also a positive semi-definite rank three matrix. Assuming that the matrix $\mathbf{MAM}$ also has at least three nonnegative eigenvalues (as is almost always the case in practice), we shall denote its three largest eigenvalues by $\lambda_1 \geq \lambda_2 \geq \lambda_3 \geq 0$, and their corresponding eigenvectors by $\mathbf{u}_1, \mathbf{u}_2, \mathbf{u}_3$. Then a classical theorem due to Sylvester, Eckart and Young on fixed-rank matrix approximations shows that the *STRAIN* is minimized by the coordinates $\mathbf{X}^* = \mathbf{M}^{-1}[\lambda_1^{1/2}\mathbf{u}_1, \lambda_2^{1/2}\mathbf{u}_2, \lambda_3^{1/2}\mathbf{u}_3]$ [10]. In chemistry, this calculation is known as the *EMBED* algorithm [5]. It should be noted that a variety of reliable iterative methods exist for finding the three largest eigenvalues and associated eigenvectors of the symmetric positive semi-definite matrix $\mathbf{MAM}$ in far less time than would be needed to fully diagonalize it.

We shall now formulate the main problem of this paper, namely embedding with a rigid substructure. Suppose that we know the 3D coordinates of a subset of *M* atoms $\mathbf{X}$, and that we wish to find 3D coordinates $\mathbf{Y}$ for the remaining *N* atoms so as to minimize the overall *STRAIN*, regarded as a function of $\mathbf{Y}$ alone. If we order the atoms appro-priately, we may partition the Gram matrix into four submatrices and write this function as

$$f(\mathbf{Y}) = \frac{1}{2} \left\| \begin{bmatrix} \mathbf{W}_A & \mathbf{W}_B \\ \mathbf{W}_B^\mathsf{T} & \mathbf{W}_C \end{bmatrix} \bullet \left( \begin{bmatrix} \mathbf{XX}^\mathsf{T} & \mathbf{XY}^\mathsf{T} \\ \mathbf{YX}^\mathsf{T} & \mathbf{YY}^\mathsf{T} \end{bmatrix} - \begin{bmatrix} \mathbf{A} & \mathbf{B} \\ \mathbf{B}^\mathsf{T} & \mathbf{C} \end{bmatrix} \right) \right\|_\mathsf{F}^2, \qquad (7)$$

where $\mathbf{A} \in \mathbf{R}^{M \times M}$, $\mathbf{B} \in \mathbf{R}^{M \times N}$, $\mathbf{C} \in \mathbf{R}^{N \times N}$, while matrices of weights satisfy $\mathbf{W}_A, \mathbf{W}_B, \mathbf{W}_C \geq \mathbf{0}$. Since $\mathbf{A} = \mathbf{XX}^\mathsf{T}$ by definition, eq. (7) may be expanded and simplified as

$$f(\mathbf{Y}) = \left\| \mathbf{W}_B \bullet [\mathbf{XY}^\mathsf{T} - \mathbf{B}] \right\|_\mathsf{F}^2 + \frac{1}{2} \left\| \mathbf{W}_C \bullet [\mathbf{YY}^\mathsf{T} - \mathbf{C}] \right\|_\mathsf{F}^2, \qquad (8)$$

which shows that $\mathbf{W}_A$ is irrelevant.

Before turning to the development of methods for finding the global minimum of the *STRAIN* function in eq. (8), we list two assumptions. First, the entries in $\mathbf{Y}$ are independent. Second, we shall consider only the unweighted problem, i.e. $w_{ij} = 1$ in both $\mathbf{W}_B$ and $\mathbf{W}_C$. It should be straightforward to extend our methods to the special case of multiplicative weights satisfying $w_{ij} = m_i m_j$, which may prove useful in some applications. The problem involving general weights will not be addressed in this paper.

# 3. First-Order Optimality Conditions

## 3.1. THE NONLINEAR MATRIX EQUATION

We begin by deriving the first-order optimality conditions for the minimization of the unweighted *STRAIN* function

$$\min_{\mathbf{Y}} (f(\mathbf{Y})) = \min_{\mathbf{Y}} \left( \|\mathbf{X}\mathbf{Y}^\mathsf{T} - \mathbf{B}\|_\mathsf{F}^2 + \tfrac{1}{2} \|\mathbf{Y}\mathbf{Y}^\mathsf{T} - \mathbf{C}\|_\mathsf{F}^2 \right). \tag{9}$$

In the following, we let $\mathbf{G_Y} = \mathbf{X}^\mathsf{T}\mathbf{X} + \mathbf{Y}^\mathsf{T}\mathbf{Y} \in \mathbf{R}^{3 \times 3}$, so that $\mathbf{G_0} = \mathbf{X}^\mathsf{T}\mathbf{X}$, and $\mathbf{V} = \mathbf{B}^\mathsf{T}\mathbf{X} \in \mathbf{R}^{N \times 3}$. Because $\mathbf{G_Y}$ is simply related to the usual inertial tensor, i.e. $tr(\mathbf{G_Y})\mathbf{I}_3 - \mathbf{G_Y}$, we shall often refer to it as *the* inertial tensor. A matrix $\mathbf{Y}_c$ at which the gradient $\nabla f(\mathbf{Y}_c) = \mathbf{0}$ shall be called a *critical matrix*.

***Theorem 1:*** A critical matrix satisfies the nonlinear matrix equation

$$\mathbf{Y}\mathbf{G_Y} - \mathbf{C}\mathbf{Y} - \mathbf{V} = \mathbf{0} \ . \tag{10}$$

Among all critical matrices, the one which maximizes the quantity

$$tr(\mathbf{Y}_c^\mathsf{T}\mathbf{V}) + \tfrac{1}{2} tr([\mathbf{Y}_c^\mathsf{T}\mathbf{Y}_c]^2) \ , \tag{11}$$

is the global minimum of $f$.

***Proof:*** Expanding eq. (9) yields

$$f(\mathbf{Y}) = tr\Big([\mathbf{XY}^\mathsf{T} - \mathbf{B}]^\mathsf{T}[\mathbf{XY}^\mathsf{T} - \mathbf{B}]\Big) + \tfrac{1}{2} tr\Big([\mathbf{YY}^\mathsf{T} - \mathbf{C}]^\mathsf{T}[\mathbf{YY}^\mathsf{T} - \mathbf{C}]\Big)$$
$$= tr\Big(\mathbf{B}^\mathsf{T}\mathbf{B} + \tfrac{1}{2}\mathbf{C}^2\Big) - tr([2\mathbf{V} + \mathbf{CY}]\mathbf{Y}^\mathsf{T}) + \qquad (12)$$
$$+ tr(\mathbf{YG_Y Y}^\mathsf{T} - \tfrac{1}{2}\mathbf{YY}^\mathsf{T}\mathbf{YY}^\mathsf{T}) .$$

The gradient $\nabla f(\mathbf{Y})$ can be arranged in the same matrix form as $\mathbf{Y}$, namely

$$\nabla f(\mathbf{Y}) \equiv \frac{\partial f}{\partial \mathbf{Y}} \equiv \begin{bmatrix} \partial f/\partial y_{11} & \partial f/\partial y_{12} & \partial f/\partial y_{13} \\ \ldots & \ldots & \ldots \\ \partial f/\partial y_{N1} & \partial f/\partial y_{N2} & \partial f/\partial y_{N3} \end{bmatrix}, \qquad (13)$$

which facilitates calculation of the gradients of matrix functions [11]. In particular, the matrix form of $\nabla tr(\mathbf{M}^\mathsf{T}\mathbf{N})$ with respect to $\mathbf{M}$ is

$$\frac{\partial tr(\mathbf{M}^\mathsf{T}\mathbf{N})}{\partial \mathbf{M}} = \frac{\partial tr(\mathbf{N}^\mathsf{T}\mathbf{M})}{\partial \mathbf{M}} = \mathbf{N} \qquad (14)$$

for any two conformable matrices $\mathbf{M}$ and $\mathbf{N}$. If we differentiate each term in eq. (12), using both eq. (14) and the rule for differentiating a matrix product [11], and then simplify using the invariance of the trace of a matrix product under cyclic permutations, we obtain eq. (10).

Given a solution $\mathbf{Y}_c$ of eq. (10), we may replace $\mathbf{Y}\mathbf{G}_\mathbf{Y}$ in the last term of eq. (12) with $\mathbf{C}\mathbf{Y}_c + \mathbf{V}$ and simplify to get

$$f(\mathbf{Y}_c) = tr\left(\mathbf{B}^\mathsf{T}\mathbf{B} + \frac{1}{2}\mathbf{C}^2\right) - tr(\mathbf{Y}_c^\mathsf{T}\mathbf{V}) - \frac{1}{2}tr([\mathbf{Y}_c^\mathsf{T}\mathbf{Y}_c]^2) \ , \tag{15}$$

thus showing why the quantity in eq. (11) must be maximized. ∎

Two comments are in order. First, because $\mathbf{G}_\mathbf{Y} = \mathbf{X}^\mathsf{T}\mathbf{X} + \mathbf{Y}^\mathsf{T}\mathbf{Y}$, eq. (10) constitutes a system of $3N$ cubic equations in the entries of $\mathbf{Y}$. The direct solution of eq. (10) is not straightforward, and alternative approaches are essential if the global minimum is to be found reliably. Second, if there are no fixed coordinates, then eq. (10) becomes

$$\mathbf{C}\mathbf{Y} = \mathbf{Y}\mathbf{G}_\mathbf{Y} . \tag{16}$$

The matrix $\mathbf{G}_\mathbf{Y} = \mathbf{Y}^\mathsf{T}\mathbf{Y}$ is a $3 \times 3$ symmetric positive-definite matrix whose spectral decomposition is denoted by $\mathbf{G}_\mathbf{Y} = \mathbf{R}\Gamma\mathbf{R}^\mathsf{T}$, $\gamma_k$ being the $k$-th diagonal entry in $\Gamma$ and $k = 1, 2, 3$. If one defines $\hat{\mathbf{Y}} = \mathbf{Y}\mathbf{R}$, so that $\hat{\mathbf{y}}_k$ is its $k$-th column, then eq. (16) becomes $\mathbf{C}\hat{\mathbf{Y}} = \hat{\mathbf{Y}}\Gamma$ or $\mathbf{C}\hat{\mathbf{y}}_k = \hat{\mathbf{y}}_k\gamma_k$.

Thus each $\hat{\mathbf{y}}_k$ is proportional to an eigenvector of $\mathbf{C}$, and since $\Gamma = \hat{\mathbf{Y}}^\top \hat{\mathbf{Y}}$ the constants of proportionality must be $\sqrt{\gamma_k}$. Finally, $f(\mathbf{Y}) = f(\hat{\mathbf{Y}})$ is minimized when $tr([\hat{\mathbf{Y}}^\top \hat{\mathbf{Y}}]^2) = \sum_k \gamma_k^2$ is maximized. This is the theorem of Sylvester, Eckart and Young mentioned in section 2.

## 3.2. EQUATIONS FOR THE INERTIAL TENSOR

In this section we derive a new system of equations and show its relation to matrix equation (10). The result is that the number of variables needed to solve minimization problem is dramatically reduced from $3N$ to $6$. We start from the spectral factorization of the $3 \times 3$ inertial tensor $\mathbf{G_Y}$,

$$\mathbf{G_Y} = \mathbf{G_0} + \mathbf{Y}^\top \mathbf{Y} = \mathbf{R}\Gamma\mathbf{R}^\top, \qquad (17)$$

order the diagonal entries of $\Gamma$ as $\gamma_1 \geq \gamma_2 \geq \gamma_3 > 0$ and, if necessary, change the signs of the columns of $\mathbf{R} = [\mathbf{r}_1, \mathbf{r}_2, \mathbf{r}_3]$ so that it is a proper orthogonal matrix. The twelve unknown entries in $\Gamma$ and $\mathbf{R}$ will be the main variables throughout the rest of this paper. It turns out that a few largest eigenvalues of a symmetric matrix $\mathbf{C}$ play a significant role so we write its spectral factorization as $\mathbf{C} = \mathbf{Q}\Sigma\mathbf{Q}^\top$ with the diagonal entries of $\Sigma$ ordered as $\sigma_1 \geq \sigma_2 \geq \ldots \geq \sigma_N$ while $\mathbf{Q}$ is an orthogonal matrix. We shall use standard notation for Kronecker delta $\delta_{kl}$, and $\mathbf{I}$ (or $\mathbf{I}_N$) for $N \times N$ identity matrix. We now define $N \times 3$ matrix $\mathbf{W} = \mathbf{Q}^\top \mathbf{V}$ whose $i$-th row is denoted by $\mathbf{w}_i^\top$, and two kinds of $3 \times 3$ variable symmetric matrices that will frequently occur, namely, a univariate one,

$$\mathbf{S}(\gamma) = \mathbf{G_0} + \mathbf{W}^\top [\gamma \mathbf{I} - \Sigma]^{-2} \mathbf{W} = \mathbf{G_0} + \sum_{i=1}^{N} \frac{\mathbf{w}_i \mathbf{w}_i^\top}{(\gamma - \sigma_i)^2} , \qquad (18)$$

and a bivariate one

$$\begin{aligned}\mathbf{S}(\gamma, \delta) &= \mathbf{G_0} + \mathbf{W}^\top ([\gamma \mathbf{I} - \Sigma]^{-1} [\delta \mathbf{I} - \Sigma]^{-1}) \mathbf{W} \\ &= \mathbf{G_0} + \sum_{i=1}^{N} \frac{\mathbf{w}_i \mathbf{w}_i^\top}{(\gamma - \sigma_i)(\delta - \sigma_i)} .\end{aligned} \qquad (19)$$

***Theorem 2:*** If $\mathbf{Y}_c$ is a critical matrix of $f$, then the eigenvalues and eigenvectors $\{\gamma_1, \mathbf{r}_1; \gamma_2, \mathbf{r}_2; \gamma_3, \mathbf{r}_3\}$ of the associated inertial tensor $\mathbf{G}_{\mathbf{Y}_c}$ satisfy the system of six equations

$$\begin{aligned}\gamma_k - \mathbf{r}_k^\top \mathbf{S}(\gamma_k) \mathbf{r}_k &= 0 , & k &= 1, 2, 3 \\ \mathbf{r}_k^\top \mathbf{S}(\gamma_k, \gamma_l) \mathbf{r}_l &= 0 , & 1 &\leq k < l \leq 3 ,\end{aligned} \qquad (20)$$

as well as the six orthonormality conditions

$$\mathbf{r}_k^\mathsf{T} \mathbf{r}_l = \delta_{kl} \qquad (21)$$

***Proof:*** The orthonormality relations (21) hold for the eigenvectors of any symmetric matrix. By premultiplying each eigenvector equation $\mathbf{G_Y} \mathbf{r}_k = \gamma_k \mathbf{r}_k$ by $\mathbf{r}_l^\mathsf{T}$, we observe that there are only six distinct such products for $1 \le k \le l \le 3$. These products are equivalent, by eqs. (17) and (21), to the following system of six equations

$$\delta_{kl} \gamma_k = \mathbf{r}_k^\mathsf{T} \mathbf{G_Y} \mathbf{r}_l = \mathbf{r}_k^\mathsf{T} \mathbf{G_0} \mathbf{r}_l + (\mathbf{Y}\mathbf{r}_k)^\mathsf{T}(\mathbf{Y}\mathbf{r}_l) \ . \qquad (22)$$

To obtain an expression for $\mathbf{Y}_c \mathbf{r}_k$ we multiply eq. (10) with the matrix $\mathbf{R}$ and substitute eq. (17) for $\mathbf{G_Y}$. The result is $\mathbf{Y}_c \mathbf{R} \Gamma = \mathbf{C} \mathbf{Y}_c \mathbf{R} + \mathbf{V} \mathbf{R}$ which is equivalent to the three vector equations $(\gamma_k \mathbf{I} - \mathbf{C}) \mathbf{Y}_c \mathbf{r}_k = \mathbf{V} \mathbf{r}_k$, one for each column of $\mathbf{R}$, whose formal solutions are given by

$$\mathbf{Y}_c \mathbf{r}_k = [\gamma_k \mathbf{I} - \mathbf{C}]^{-1} \mathbf{V} \mathbf{r}_k \qquad k = 1, 2, 3 . \qquad (23)$$

Substitution of eq. (23) into eq. (22) and rearrangement results in

$$\delta_{kl} \gamma_k - \mathbf{r}_k^\mathsf{T} [\mathbf{G_0} + \mathbf{V}^\mathsf{T}([\gamma_k \mathbf{I} - \mathbf{C}]^{-1}[\gamma_l \mathbf{I} - \mathbf{C}]^{-1}) \mathbf{V}] \mathbf{r}_l = 0 \ . \qquad (24)$$

Separation of these equations into the $k = l$ and $k < l$ cases, and substitution of $\mathbf{C} = \mathbf{Q}\Sigma\mathbf{Q}^\top$ and $\mathbf{V} = \mathbf{Q}\mathbf{W}$, now yields eqs. (20). ∎

The equations (20) and (21) can be regarded as a system of twelve equations in the twelve unknowns $\{\gamma_1, \mathbf{r}_1; \gamma_2, \mathbf{r}_2; \gamma_3, \mathbf{r}_3\}$. The six constraints in eq. (21) on the nine entries of $\mathbf{R}$ can be eliminated by a suitable parametrization for the $3 \times 3$ orthogonal matrices. Two specific parametrizations will be introduced later on, the outcomes being two different nonlinear systems of six equations in six unknowns. We shall often refer to the equations (20) as the *inertial equations,* and denote them by $g_{kl}$, for $1 \le k \le l \le 3$. An eigenvalue-eigenvector pair $\{\gamma_k, \mathbf{r}_k\}$ will be referred to as an *inertial pair*. Each of the first three equations depends on a single inertial pair and will be referred to as *quadratic,* while the last three equations, each depending on two inertial pairs, will be referred to as *bilinear* ones.

The connection between inertial pairs satisfying eqs. (20)-(21) and critical matrices $\mathbf{Y}_c$ is that given any positive diagonal matrix $\Gamma$, such that each diagonal entry satisfies $0 < \gamma_k \ne \sigma_i$, and an arbitrary $3 \times 3$ orthogonal matrix $\mathbf{R}$, one can compute the matrix $\mathbf{R}\Gamma\mathbf{R}^\top$ and the $N \times 3$ matrix

$$\mathbf{Y} = [\ (\gamma_1 \mathbf{I} - \mathbf{C})^{-1} \mathbf{V} \mathbf{r}_1,\ (\gamma_2 \mathbf{I} - \mathbf{C})^{-1} \mathbf{V} \mathbf{r}_2,\ (\gamma_3 \mathbf{I} - \mathbf{C})^{-1} \mathbf{V} \mathbf{r}_3\ ] \mathbf{R}^T \qquad (25)$$

satisfying eq. (10). In general, however, $\mathbf{R} \Gamma \mathbf{R}^T \neq \mathbf{G}_0 + \mathbf{Y}^T \mathbf{Y}$ unless the matrix pair $\{\Gamma, \mathbf{R}\}$ is also a solution to eqs. (20), in which case $\mathbf{Y}$ is a critical matrix of $f$. The same critical matrix is obtained for each of the eight possible choices of signs $\pm \mathbf{r}_k$ ($k = 1, 2, 3$). The choice of the particular orientation will be described in sect. **5.2**. In the following, after computing solutions to eqs. (20)-(21), the corresponding critical matrix will be constructed according to eq. (25).

# 4. The Structure of the Solution Space

## 4.1. THE INTERSECTION CURVE OF A SPHERE AND AN ELLIPSOID

Let $g(\gamma, \mathbf{r}) = 0$ denote any one of the three quadratic inertial equations in (20), with the associated matrix $\mathbf{S}(\gamma)$ defined by eq. (18), and let $\{\gamma, \mathbf{r}\}$ denote an inertial pair that satisfies it. A geometric characterization of all such inertial pairs comes from the following result:

***Lemma 3:*** If an inertial pair $\{\gamma, \mathbf{r}\}$ satisfies the quadratic inertial equation $g(\gamma, \mathbf{r}) = 0$, then $\mathbf{r}$ is located on the intersection curve between the unit sphere and the ellipsoid defined by the matrix $\mathbf{S}(\gamma)/\gamma$.

***Proof:*** Since $\gamma > 0$ the equation $g(\gamma, \mathbf{r}) = 0$ can be rewritten as

$$\mathbf{r}^\mathsf{T} \left[ \frac{\mathbf{S}(\gamma)}{\gamma} \right] \mathbf{r} = 1 \ . \tag{26}$$

For any fixed $0 < \gamma \neq \sigma_i$, the matrix $\mathbf{S}(\gamma)/\gamma$ is positive-definite, and any vector $\mathbf{r}$ that satisfies eq. (26) lies on the surface of an ellipsoid centered at the origin. As a column of an orthogonal matrix, however, the vector $\mathbf{r}$ also lies on the surface of the unit sphere centered at origin. Therefore, $\mathbf{r}$ lies on the intersection curve between the unit sphere and the ellipsoid defined by $\gamma$ via eq. (26) (see Fig. 1). ∎

It is well known that the spectral factorization of $\mathbf{S}(\gamma)/\gamma = \mathbf{U}\Lambda\mathbf{U}^\mathsf{T}$ determines the principal axes of this ellipsoid as the three eigenvectors in $\mathbf{U}$, while the three eigenvalues $0 < \lambda_3 \leq \lambda_2 \leq \lambda_1$ are the squares of the three semi-axes. It is also evident from eqs. (18) and (26) that both $\Lambda = \Lambda(\gamma)$ and $\mathbf{U} = \mathbf{U}(\gamma)$ are continuous functions, except at the eigenvalues of $\mathbf{C}$ (see also [12]). Therefore, the problem of existence of an intersection curve, henceforth denoted by $\mathbf{r}(\gamma)$, consists of finding inter-vals of γ over which the solutions of the quadratic vector equation (26) constitute real curves.

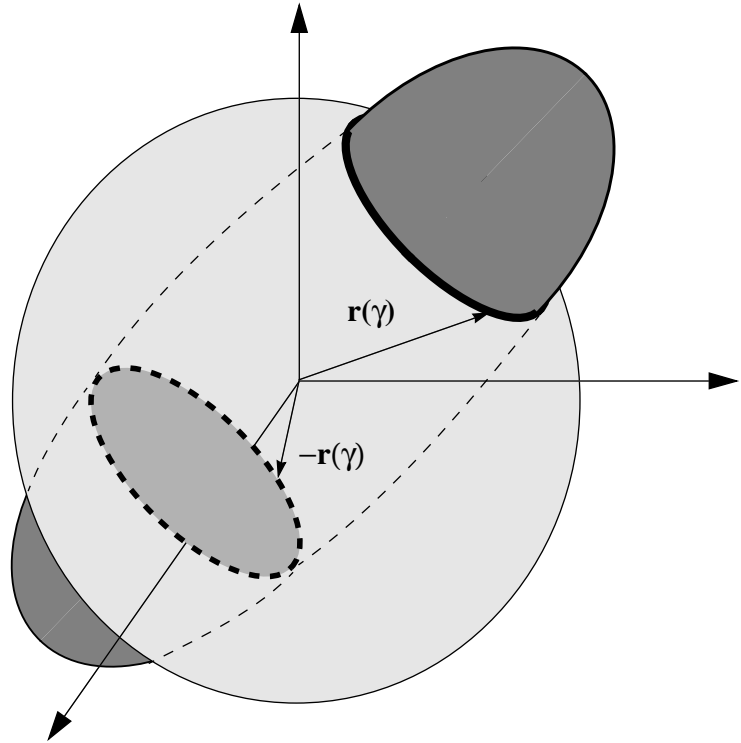

**Figure 1. The two symmetric curves in which the ellipsoid defined by the matrix $S(\gamma)/\gamma$ intersects the unit sphere.**

Since the ellipsoid in eq. (26) and the unit sphere are centrally symmetric convex surfaces, the necessary and sufficient condition for their intersection is

$$\lambda_3(\gamma) \leq 1 \leq \lambda_1(\gamma) \ . \tag{27}$$

In other words, there is no intersection if all semi-axes are strictly greater, or strictly smaller, than the radius of the sphere. By the central symmetry of both surfaces, the intersection consists of a symmetrical pair of (nonintersecting) closed spherical curves, related to each other by reflection in the origin. The curves are symmetric with respect to all three axes, but we shall be mostly interested in the axis of symmetry along which piercing of surfaces occurs: $\mathbf{u}_1$ if $\lambda_2(\gamma) < 1$, or $\mathbf{u}_3$ if $1 < \lambda_2(\gamma)$. We shall label this axis vector by $\mathbf{u}_k$, and the other two axes by $\mathbf{u}_i$ and $\mathbf{u}_j$, so that $k = 1$ implies $(i, j) = (2, 3)$ while $k = 3$ implies $(i, j) = (1, 2)$. At a *bifurcation point* $\gamma_b$, which is defined by the nonlinear equation

$$\lambda_2(\gamma_b) = 1 , \tag{28}$$

the two intersection curves merge and the axis of symmetry for both curves undergoes a discontinuous change from $\mathbf{u}_1$ to $\mathbf{u}_3$.

The following lemma shows how the explicit analytic parametrization of the intersection curve depends on the axis of symmetry.

***Lemma 4:*** The projection of both intersection curves onto the central plane orthogonal to the axis of symmetry $\mathbf{u}_k$ is an ellipse, whose semi-axes are given by

$$\begin{aligned} \beta_i(\gamma) &= \sqrt{(\lambda_k - 1)/(\lambda_k - \lambda_i)} \\ \beta_j(\gamma) &= \sqrt{(\lambda_k - 1)/(\lambda_k - \lambda_j)} \end{aligned} . \qquad (29)$$

Any regular parametrization of an ellipse can be chosen to parametrize the intersection curve. We shall use the trigonometric one given by

$$\begin{aligned} \xi_i(\gamma, \psi) &= \beta_i(\gamma) \cos(\psi) \\ \xi_j(\gamma, \psi) &= \beta_j(\gamma) \sin(\psi) \\ \xi_k(\gamma, \psi) &= \pm\sqrt{1 - \xi_i^2(\gamma, \psi) - \xi_j^2(\gamma, \psi)} \end{aligned} , \qquad (30)$$

for $\psi \in [0, 2\pi)$. The vector $\xi(\gamma, \psi) = [\xi_1(\gamma, \psi), \xi_2(\gamma, \psi), \xi_3(\gamma, \psi)]^\top$ is the representation of intersection curve with respect to the principal axes, while its representation with respect to the standard axes is given by

$$\mathbf{r}(\gamma, \psi) = \mathbf{U}(\gamma)\, \xi(\gamma, \psi) . \qquad (31)$$

***Proof:*** If we write the equations for the unit sphere and the ellipsoid in the principal axis coordinates $\xi$ of the ellipsoid, and then eliminate the coordinate $\xi_k$ from these equations, we get

$$\left(\frac{\lambda_k - \lambda_i}{\lambda_k - 1}\right)\xi_i^2 + \left(\frac{\lambda_k - \lambda_j}{\lambda_k - 1}\right)\xi_j^2 = 1 \,. \tag{32}$$

Equation (32) describes an ellipse in the plane spanned by $i$-th and $j$-th eigenvector, whose associated semi-axes are given by eq. (29). If we choose an angle $\psi$ for the trigonometric parametrization of the ellipse, and back substitute for the eliminated coordinate $\xi_k$, we get eq. (30). The vector $\xi(\gamma, \psi)$ refers to the principal axes; transformation to the original coordinates comes from the identification $\mathbf{U}^\mathsf{T}\mathbf{r} = \xi(\gamma, \psi)$, which proves eq. (31). ∎

*Lemma 4* provides a simple and efficient way to compute any point on any intersection curve where, according to *Lemma 3*, all solutions to the inertial equations are to be found. Since we shall also require the derivatives of $\mathbf{r}(\gamma, \psi)$ with respect to both variables during the search for the global minimum, we should analyze the regions where these derivatives are ill-behaving in order to prevent difficulties or even failure of the search. Geometrically this happens if a solution is located on a set of very small, or very narrow, intersection curves. Analytically, the small intersection curves occur when both projected semi-axes $\beta_i(\gamma), \beta_j(\gamma) \to 0$ simultaneously, which happens when either $\lambda_1(\gamma) \to 1$ or $\lambda_3(\gamma) \to 1$. Narrow intersection curves occur when the projected eccentricity $min((\beta_i/\beta_j), (\beta_j/\beta_i)) \ll 1$. Almost all such cases occur in the vicinity of those bifurcation points that have a very small angle between intersection curves (see Fig. 2).

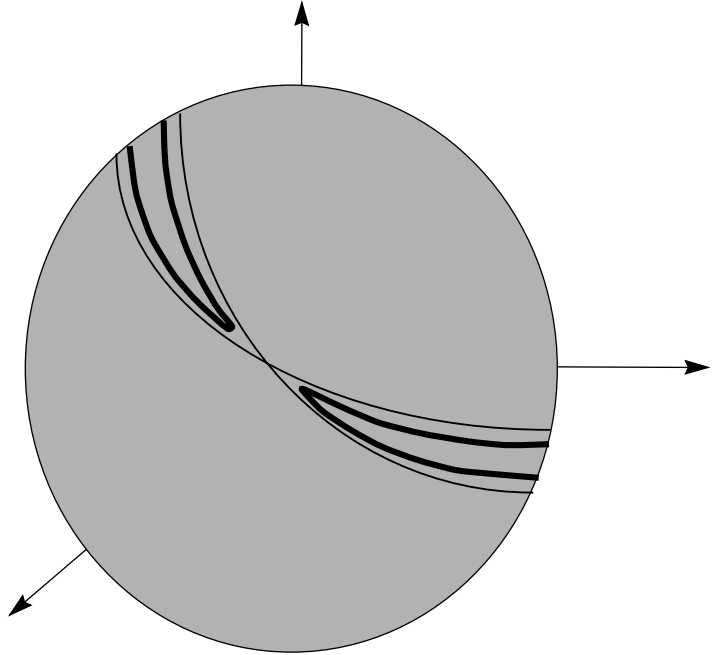

**Figure 2. The shape of intersection curves due to the narrow ellipsoids near a bifurcation point.**

To compute this angle, we eliminate $\xi_2$ from the equations for an ellipsoid and the unit sphere, set $\lambda_2(\gamma_b) = 1$, and obtain

$$tan\left(\frac{\alpha}{2}\right) = \frac{\xi_1}{\xi_3} = \pm\sqrt{\frac{1-\lambda_3(\gamma_b)}{\lambda_1(\gamma_b)-1}} . \tag{33}$$

This shows that the very small angles are caused by either $\lambda_3(\gamma_b)$ being very close to 1 (which is also a near-tangency case), or by $\lambda_1(\gamma_b) \gg 1$ (which happens when a bifurcation point is very close to an eigenvalue $\sigma_k$), or both.

## 4.2. SPHERICAL BOUNDS ON THE MOMENTS OF INERTIA

Let $\{\gamma, \mathbf{r}\}$ be an unknown inertial pair satisfying the quadratic equation $g(\gamma, \mathbf{r}) = 0$. Since $\mathbf{r}$ is confined to the unit sphere, it is of great interest to find an upper bound on its eigenvalue $\gamma > 0$ as well as a better lower bound, because such bounds would reduce the size of the subspace over which the search for the global minimum takes place. The notation $\lambda_{max}(\cdot)$ and $\lambda_{min}(\cdot)$ stands for the maximal and minimal eigenvalue of the argument matrix (e.g. from eq. (26), $\lambda_{max}(\mathbf{S}(\gamma)) = \gamma\lambda_1(\gamma)$).

**Theorem 5:** Let $\gamma^+$ be the maximal solution of the nonlinear equation

$$\gamma = \lambda_{max}(\mathbf{S}(\gamma)) \ . \tag{34}$$

Then the largest eigenvalue $\gamma_1$ of the inertial tensor is bounded by

$$max(\sigma_1, \lambda_{max}(\mathbf{G_0})) \leq \gamma_1 \leq \gamma^+ . \tag{35}$$

***Proof:*** We start by writing the quadratic inertial equation (26) as $\gamma = \mathbf{r}^\mathsf{T}\mathbf{S}(\gamma)\mathbf{r}$. The maximum of its right-hand side over all unit vectors $\mathbf{r}$, (also known as the Rayleigh quotient of $\mathbf{S}(\gamma)$) is $\lambda_{max}(\mathbf{S}(\gamma))$ (see e.g. [13]), and the resulting inequality

$$\gamma \leq \max_{\|\mathbf{r}\|=1}\left(\mathbf{r}^\mathsf{T}\mathbf{S}(\gamma)\mathbf{r}\right) = \lambda_{max}(\mathbf{S}(\gamma)) , \tag{36}$$

yields the nonlinear equation (34). If the maximal positive solution of eq. (34) exists, then it provides an upper bound on any eigenvalue $\gamma$, and in particular on the largest one $\gamma_1$.

To prove the existence of a solution to eq. (34), we recall that $\mathbf{S}(\gamma)$ is a positive-definite matrix, so that maximal eigenvalue function

$$\rho(\gamma) = \lambda_{max}(\mathbf{S}(\gamma)) = \lambda_{max}\left(\mathbf{G_0} + \sum_{i=1}^{N} \frac{\mathbf{w}_i\mathbf{w}_i^\mathsf{T}}{(\gamma-\sigma_i)^2}\right) \tag{37}$$

is strictly positive. It is also a rational function with poles of order two at the eigenvalues of $\mathbf{C}$. In the interval $(\sigma_1, \infty)$, it is convex and decreases monotonically from infinity at the pole $\sigma_1$ to the asymptotic value of $\lambda_{max}(\mathbf{G_0})$. We omit the proofs of these statements, which are based on the negative and positive semi-definiteness of the first and second derivatives of $\mathbf{S}(\gamma)$, respectively. It follows that in the interval $(\sigma_1, \infty)$ there is only one solution, denoted by $\gamma^+$, to the equation $\gamma = \rho(\gamma)$, which is also the maximal solution to the inequality $\gamma \leq \rho(\gamma)$.

To prove a lower bound on $\gamma_1$ in eq. (35), consider first the case when $\sigma_1 \leq \lambda_{max}(\mathbf{G_0})$. Then, rewriting eq. (17) as $\Gamma = \mathbf{R}^\mathsf{T}(\mathbf{G_0} + \mathbf{Y}^\mathsf{T}\mathbf{Y})\mathbf{R}$, we see that $\gamma_1 = \lambda_{max}((\mathbf{R}^\mathsf{T}(\mathbf{G_0} + \mathbf{Y}^\mathsf{T}\mathbf{Y})\mathbf{R})) = \lambda_{max}(\mathbf{G_0} + \mathbf{Y}^\mathsf{T}\mathbf{Y}) \geq \lambda_{max}(\mathbf{G_0})$, by the positive semi-definiteness of $\mathbf{Y}^\mathsf{T}\mathbf{Y}$ and the positive definiteness of $\mathbf{G_0}$. In the case when $\lambda_{max}(\mathbf{G_0}) < \sigma_1$, let $\{\gamma_1, \mathbf{r}_1\}$ be the first inertial pair satisfying the quadratic eq. (20) and proceed as follows:

$$\begin{aligned}\gamma_1 &= \mathbf{r}_1^\mathsf{T}\left(\mathbf{G_0} + \sum_{i=1}^{N} \frac{\mathbf{w}_i \mathbf{w}_i^\mathsf{T}}{(\gamma_1 - \sigma_i)^2}\right)\mathbf{r}_1 \geq \mathbf{r}_1^\mathsf{T}\mathbf{G_0}\mathbf{r}_1 + \frac{(\mathbf{r}_1^\mathsf{T}\mathbf{w}_1)^2}{(\gamma_1 - \sigma_1)^2} \\ &\geq \lambda_{min}(\mathbf{G_0}) + \frac{(\mathbf{r}_1^\mathsf{T}\mathbf{w}_1)^2}{(\gamma_1 - \sigma_1)^2}\end{aligned} \qquad (38)$$

The first inequality comes from keeping only the first term under the sum (all terms are nonnegative), while the second one comes from minimizing the Rayleigh quotient of $\mathbf{G_0}$. By rewriting the final inequality in (38) as $(\gamma_1 - \sigma_1)^3 + (\sigma_1 - \lambda_{min}(\mathbf{G_0}))(\gamma - \sigma_1)^2 \geq (\mathbf{r}_1^\mathsf{T}\mathbf{w}_1)^2$ and defining

$\delta = \gamma_1 - \sigma_1$, we obtain the cubic inequality

$$\delta^3 + (\sigma_1 - \lambda_{min}(\mathbf{G_0}))\delta^2 - (\mathbf{r}_1^T\mathbf{w}_1)^2 \geq 0. \tag{39}$$

There is always one nonnegative solution to the inequality (39) (see discussion following eq. (42)), which means that $\gamma_1 \geq \sigma_1$. This would complete the proof for the lower bound in eq. (35), were it not for occasional occurrence of two negative real solutions to the inequality (39). We claim that they can be ignored since only the greatest of the available lower bounds is of interest. ∎

The inequality (35) says that there is no intersection between the variable ellipsoid and the unit sphere as long as $\gamma > \lambda_{max}(\mathbf{S}(\gamma))$. The bound $\gamma^+$ is the tangency point at which intersection commences for the first time, and so we have named it the "spherical bound".

We shall solve equation (34) by a rational interpolation method which consists of two iterations. An outer iteration calculates the coefficients of a rational interpolant to function $\rho(\gamma)$, defined by eq. (37), and monitors convergence to the solution of eq. (34), while an inner one computes the maximal root of a special cubic equation similar to (39).

The outer iteration interpolates $\rho(\gamma)$ by the simplest rational function having a pole of order two at $\sigma_1$, namely

$$\hat{\rho}(\gamma) = c + \frac{d}{(\gamma - \sigma_1)^2} . \qquad (40)$$

The coefficients $c, d$ are computed from the interpolatory data $\rho^{i-1} = \rho(\gamma^{i-1})$ and $\rho^i = \rho(\gamma^i)$ at the last two iterates $\gamma^{i-1}$ and $\gamma^i$, as

$$d = \frac{\rho^i - \rho^{i-1}}{(\gamma^i - \sigma_1)^{-2} - (\gamma^{i-1} - \sigma_1)^{-2}} , \qquad c = \rho^i - \frac{d}{(\gamma^i - \sigma_1)^2} . \qquad (41)$$

The next iterate $\gamma^{i+1}$ is obtained by solving the equation $\gamma = \hat{\rho}(\gamma)$ to working precision, which is the task of the inner iteration. Given the current values of the coefficients $c, d$, and with the change of variable $\delta = \gamma - \sigma_1$, eq. (34) takes on the following form:

$$\delta^3 + (\sigma_1 - c)\delta^2 - d = 0 . \qquad (42)$$

Since $\rho(\gamma)$ is a strictly decreasing function for $\gamma > \sigma_1$, it follows from eq. (41) that $d > 0$ for all iterates $\gamma^{i-1}, \gamma^i > \sigma_1$. Therefore, the cubic polynomial in eq. (42) is strictly negative at $\delta = 0$, which guarantees that it has a single positive root. A simple transformation of eq. (42) yields the following bounds on this root:

$$\sqrt{\frac{d}{\sigma_1 - c + \sqrt{\frac{d}{\sigma_1 - c}}}} < \delta < \sqrt{\frac{d}{\sigma_1 - c}}, \quad \text{if } \sigma_1 - c > 0$$

$$|\sigma_1 - c| < \delta < |\sigma_1 - c| + \frac{d}{(\sigma_1 - c)^2}, \quad \text{if } \sigma_1 - c < 0.$$

(43)

For any point in the interval (43), both the first and the second derivative of the cubic polynomial in eq. (42) are positive. This means that Newton's method converges unconditionally and monotonically, except possibly for the first step, to the positive root $\delta^i$ of eq. (42). This completes the inner iteration, after which we set $\gamma^{i+1} = \sigma_1 + \delta^i$ in the outer one.

By using the convexity of $\rho(\gamma)$ in the interval $(\sigma_1, \infty)$, it can be shown that this rational interpolation procedure generates a sequence $\gamma^i$ which converges superlinearly to a solution of $\rho(\gamma) = \gamma$ starting from any pair of estimates $\gamma^1, \gamma^2 > \sigma_1$ [15]. In practice, the convergence rate tends to be nearly quadratic since $\hat{\rho}(\gamma)$ inherits the convexity of $\rho(\gamma)$.

A good starting point for the outer iteration is obtained by computing an upper bound on the right-hand side of eq. (37) as

$$\rho(\gamma) \leq \lambda_{max}(\mathbf{G_0}) + \lambda_{max}\left(\sum_{i=1}^{N} \frac{\mathbf{w}_i \mathbf{w}_i^\mathsf{T}}{(\gamma - \sigma_i)^2}\right)$$
$$\leq \lambda_{max}(\mathbf{G_0}) + \sum_{i=1}^{N} \frac{\mathbf{w}_i^\mathsf{T} \mathbf{w}_i}{(\gamma - \sigma_i)^2} \leq \lambda_{max}(\mathbf{G_0}) + \frac{\|\mathbf{W}\|_F^2}{(\gamma - \sigma_1)^2} \ . \quad (44)$$

The last inequality is valid for all $\gamma \in (\sigma_1, \infty)$. On cross-multiplying and setting $\delta = \gamma - \sigma_1$, we obtain a cubic inequality

$$\delta^3 + (\sigma_1 - \lambda_{max}(\mathbf{G_0}))\delta^2 - \|\mathbf{W}\|_F^2 \leq 0 \ , \quad (45)$$

quite similar to eqs. (42) and (39). The maximum real solution $\delta_0$ to eq. (45) is strictly positive, so that $\gamma_1 \leq \gamma^+ \leq \sigma_1 + \delta_0$. The upper bound $\sigma_1 + \delta_0$ is generally already quite good and, together with any other point in the interval $(\sigma_1, \sigma_1 + \delta_0)$, provides the desired starting pair for the iteration in eqs. (40)-(41).

A lower bound $\bar{\gamma}$ on all the eigenvalues $\gamma$ is obtained by minimizing the Rayleigh quotient in $\gamma = \mathbf{r}^\mathsf{T} \mathbf{S}(\gamma) \mathbf{r}$. The result is

$$\gamma \geq \lambda_{min}(\mathbf{S}(\gamma)) = \mu(\gamma) \quad (46)$$

We use a rational interpolation procedure, similar to that in eqs. (40)-(41) to solve for the minimum value $\gamma^-$ such that $\mu(\gamma^-) = \gamma^-$. Knowing that $\gamma^- \geq \lambda_{min}(\mathbf{G}_0)$, we start the search from $\lambda_{min}(\mathbf{G}_0)$. The meaning of eq. (46) is that the ellipsoid $\mathbf{S}(\gamma)/\gamma$ is guaranteed to have all its semi-axes less than 1 for any $\gamma < \gamma^-$, so there are no further intersections with the unit sphere.

## 4.3. CONICAL BOUNDS ON THE MOMENTS OF INERTIA

By definition, there are two antipodal tangency points at $\gamma = \gamma^+$, where the ellipsoid touches the unit sphere contained within it. As $\gamma$ starts decreasing from $\gamma^+$ towards $\sigma_1$, so does the size of the ellipsoid, because all three semi-axes decrease. The two intersection curves start growing until either the pole $\sigma_1$ or a bifurcation point $\gamma_b$ is reached, if there is one in this interval. At the bifurcation point the intersection curves meet at a single pair of antipodal points. On further decreasing the value of $\gamma$, the curves bifurcate (switch to complementary regions defined by eq. (33), see Fig. 2), but now begin to shrink in size. An important property of the continuous variation of the intersection curves $\mathbf{r}(\gamma)$ before and after $\gamma_b$ is that they are *nested,* as shown in Fig. 3.

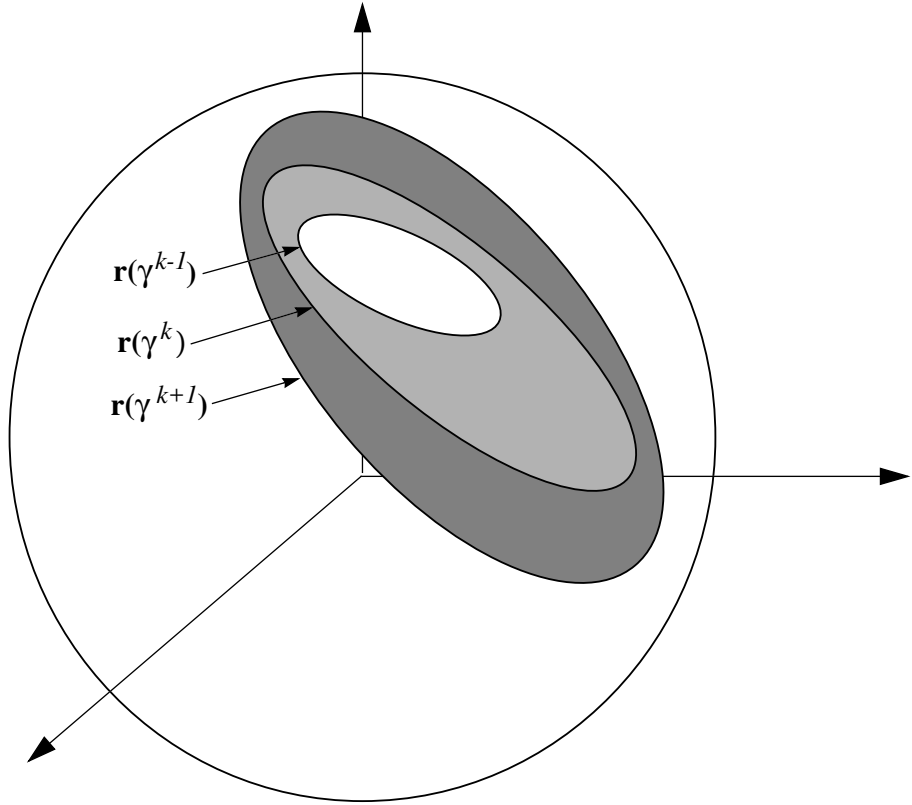

**Figure 3. A nested sequence of intersection curves defined by a decreasing sequence of $\gamma$ values.**

***Lemma 6:*** If there is no bifurcation point in the interval $(\sigma_1, \gamma^+)$, then for any $\sigma_1 < \gamma'' < \gamma' < \gamma^+$ the two intersection curves $\mathbf{r}(\gamma')$ and $\mathbf{r}(\gamma'')$ are nested on the spherical surface, in the sense that (a) they do not intersect, and (b) the set of points on the unit sphere bounded by $\mathbf{r}(\gamma'')$ contains the set of points bounded by $\mathbf{r}(\gamma')$, including the initial tangency point $\mathbf{r}(\gamma^+)$. If, on the other hand, there is a bifurcation point in this interval, i.e. $\sigma_1 < \gamma_b$, then the statement above continues to hold for $\gamma_b < \gamma'' < \gamma' < \gamma^+$, while in the domain $\sigma_1 < \gamma'' < \gamma' < \gamma_b$ the containment is reversed: the set of points on the unit sphere bounded by $\mathbf{r}(\gamma'')$ is contained in the set of points bounded by any previous $\mathbf{r}(\gamma')$ including the limit curve $\mathbf{r}(\gamma_b)$.

***Proof:*** Starting with two distinct intersection curves at $\gamma' \neq \gamma''$, one assumes that they intersect at $\mathbf{p} = \mathbf{r}(\gamma', \psi') = \mathbf{r}(\gamma'', \psi'')$. By subtracting the quadratic equation $\mathbf{p}^T \mathbf{S}(\gamma')\mathbf{p} = \gamma'$ from $\mathbf{p}^T \mathbf{S}(\gamma'')\mathbf{p} = \gamma''$ and simplifying, one gets the expression

$$1 = \sum_{k=1}^{N} (\mathbf{w}_k^T \mathbf{p})^2 \left( \frac{2\sigma_k - \gamma' - \gamma''}{(\gamma' - \sigma_k)^2 (\gamma'' - \sigma_k)^2} \right) \quad (47)$$

whose right hand side is always negative for any choice of $\sigma_1 < \gamma', \gamma''$. The contradiction in (47) can be traced back to the assumption that the curves intersect, which completes the proof on the nesting property. The two cases of containment are direct consequences of the inequalities among the projected semi-axes in eq. (29), before and after the bifurcation point. ∎

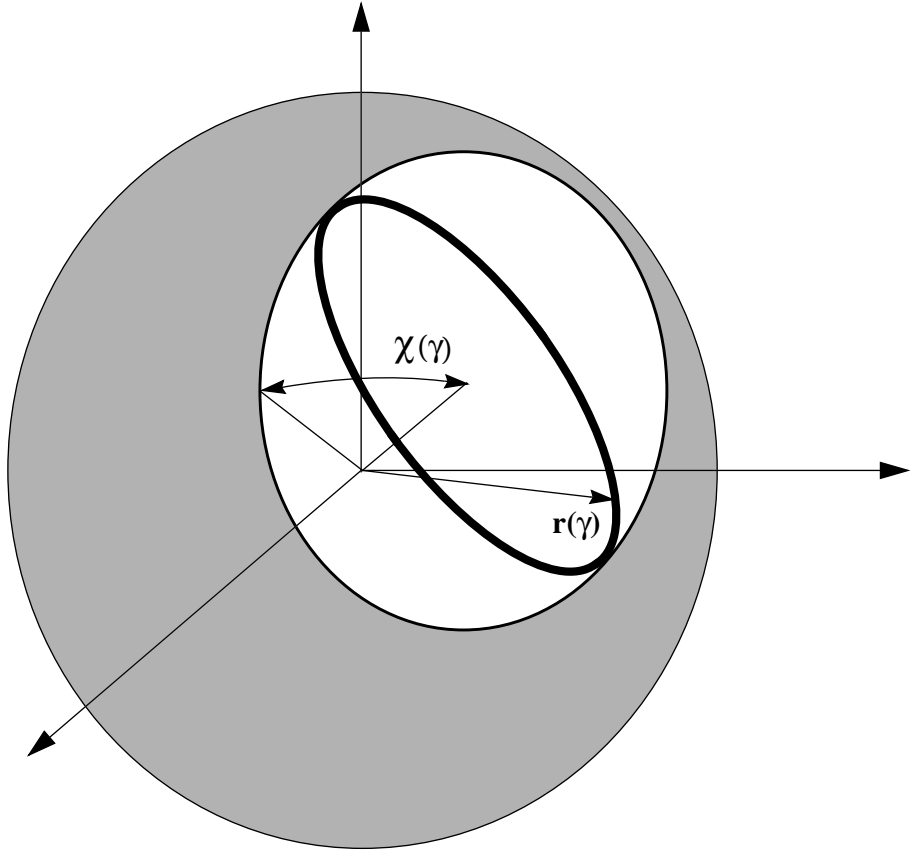

**Figure 4. The cone defined by the intersection curve r(γ) and its half-angle χ(γ).**

Let $\chi(\gamma)$ be the maximum angle between any vector on the intersection curve $\mathbf{r}(\gamma)$ and its axis of symmetry. In view of eq. (29), this angle is given by

$$\chi(\gamma) = sin^{-1}(max(\beta_i, \beta_j)) \ . \tag{48}$$

We shall now determine an upper bound $\gamma_2'$ on the second eigenvalue by analyzing this angle. To start with, it can be identified with the half-angle of the circular cone having the same symmetry axis as $\mathbf{r}(\gamma)$ (see Fig. 4).

In other words, the angle between any two points on the intersection curve $\mathbf{r}(\gamma)$ does not exceed $2\chi(\gamma)$. For that reason, we named $\gamma_2^+$ the "conical" bound. In what follows, we shall need the following definition:

$$\gamma_{b_1} = \begin{cases} \gamma_b & \text{if } \sigma_1 < \gamma_b \\ \sigma_1 & \text{otherwise} \end{cases}. \tag{49}$$

***Theorem 7:*** The second eigenvalue is bounded from above by the unique solution $\gamma_2^+$ of the equation

$$\lambda_1(\gamma) + \lambda_2(\gamma) = 2 \tag{50}$$

in the interval $(\gamma_{b_1}, \gamma^+)$.

***Proof:*** By taking the sine of both sides of eq. (48), followed by squaring, making use of eq. (29) (recall that the eigenvalues of the matrix $\mathbf{S}(\gamma)/\gamma$ are ordered as $\lambda_3(\gamma) \leq \lambda_2(\gamma) \leq \lambda_1(\gamma)$) and rearranging, we find that the condition $\chi(\gamma) = \pi/4$ is equivalent to eq. (50).

We now show that equation $\chi(\gamma) = \pi/4$ is a consequence of the orthogonality among eigenvectors. By **Lemma 6**, the function $\chi(\gamma)$ is monotonically decreasing such that at the tangency point $\chi(\gamma^+) = 0$, while at a bifurcation point $\chi(\gamma_b) = \pi/2$, which follows from eq. (55) and the fact that $max(\beta_2, \beta_3) = 1$ at $\lambda_2 = 1$ by eq. (29). If there is no bifurcation point in the interval $(\sigma_1, \gamma^+)$, then $\sigma_1$ acts like one because when $\gamma \to \sigma_1$ then $\beta_2(\gamma) \to 1$ and $\chi(\gamma) \to \pi/2$. Hence, there is a unique value $\gamma_2^+ \in (\gamma_{b_1}, \gamma^+)$ at which $\chi(\gamma_2^+) = \pi/4$.

To show that $\gamma_2^+$ is an upper bound on $\gamma_2$, suppose first that the maximal eigenvalue satisfies $\gamma_1 \leq \gamma_2^+$, which does happen albeit infrequently. In that case we must also have $\gamma_2 \leq \gamma_2^+$, since $\gamma_2 \leq \gamma_1$ by the ordering of the eigenvalues. In the more frequent case when $\gamma_2^+ < \gamma_1$, assume that $\gamma_2^+ < \gamma_2 \leq \gamma_1$. Then by **Lemma 6**, the set of points within the intersection curve $\mathbf{r}(\gamma_1)$ is contained in the set of points within the intersection curve $\mathbf{r}(\gamma_2)$, known to satisfy $2\chi(\gamma_2) < \pi/2$. It follows that there is no vector on the intersection curve $\mathbf{r}(\gamma_1)$ and no vector on the intersection curve $\mathbf{r}(\gamma_2)$ that could be mutually orthogonal and thus satisfy orthogonality condition (21). This contradiction is due to the assumption that $\gamma_2^+ < \gamma_2$ and we conclude that $\gamma_2^+$ must be an upper bound on $\gamma_2$, as well as on $\gamma_3$, by the ordering of the eigenvalues. ∎

To solve the nonlinear equation (50) we use the discretization of $\gamma_1$ over the interval $(max(\sigma_1, \lambda_{max}(\mathbf{G_0})), \gamma_1^+)$ that will be described in section **6.3**. Upon bracketing the solution, we refine an initial solution by means of Newton's method. The calculation of the requisite eigenvalue derivatives $\lambda_1'$ and $\lambda_2'$ is described later on (eq. (71)). This iteration is guaranteed to converge rapidly due to the monotonicity and differentiability of the underlying function in this interval.

In summary, we have shown that $\gamma_1 \in (max(\sigma_1, \lambda_{max}(\mathbf{G_0})), \gamma^+)$ and $\gamma_2, \gamma_3 \in (\gamma^-, \gamma_2^+)$ (see eq. (46)). We note that we have never yet found an example in which $\gamma_3 < \sigma_3$ at the global minimum of $f$. This leads to the following:

***Conjecture 8:*** At the global minimum of $f$ the smallest eigenvalue satisfies $\gamma_3 > max(\sigma_3, \gamma^-)$.

The conjectural part is the case when $\gamma^- < \sigma_3$, because the other case follows easily from the definition of $\gamma^-$ in eq. (46). If this conjecture turns out to be true, we will be able to further reduce the range of values accessible to the moments of inertia and hence reduce the time needed to find the global minimum. In the absence of a proof, we make sure that we find the global minimum by starting the search from $\bar{\gamma}$, and accept some degradation of performance due to this less stringent lower bound.

# 5. Iterative Methods

## 5.1. THE NEWTON-KANTOROVITCH METHOD FOR ALL VARIABLES

The first method of solving $\Phi(\mathbf{Y}) = \nabla f(\mathbf{Y}) = \mathbf{0}$ is based on eq. (10). It is a matrix version of the Newton-Kantorovitch method for solving operator equations [14]. The method exhibits quadratic convergence, and we have observed very few examples in which it fails to converge, although it may of course converge to a saddle point, local minimum or maximum.

Given any $\mathbf{Y}_n \in \mathbf{R}^{N \times 3}$ and a perturbation matrix $\Delta \in \mathbf{R}^{N \times 3}$, we expand the function $\Phi(\mathbf{Y}_n + \Delta)$ as

$$\Phi(\mathbf{Y}_n + \Delta) = [\mathbf{Y}_n + \Delta][(\mathbf{Y}_n + \Delta)^\mathsf{T}(\mathbf{Y}_n + \Delta) + \mathbf{X}^\mathsf{T}\mathbf{X}] - \mathbf{C}[\mathbf{Y}_n + \Delta] - \mathbf{V}$$
$$= \Phi(\mathbf{Y}_n) + \mathbf{D}_\Delta(\Phi) + \mathbf{O}(\|\Delta\|^2) \, , \quad (51)$$

where $\mathbf{O}(\|\Delta\|^2)$ denotes matrices whose norm depends at least quadratically on $\|\Delta\|$. The terms linear in $\Delta$ are collected in the matrix

$$\mathbf{D}_\Delta(\Phi) = \Delta[\mathbf{Y}_n^\mathsf{T}\mathbf{Y}_n + \mathbf{X}^\mathsf{T}\mathbf{X}] + \mathbf{Y}_n \Delta^\mathsf{T} \mathbf{Y}_n + [\mathbf{Y}_n \mathbf{Y}_n^\mathsf{T} - \mathbf{C}]\Delta \, , \quad (52)$$

which is the directional or Frechet derivative of the matrix-valued function $\Phi$ evaluated in the direction $\Delta$ [14]. Setting the right-hand side of eq. (51) to the zero matrix and ignoring the higher-order terms, we obtain a system of $3N$ linear equations in the $3N$ unknown entries of $\Delta$, namely

$$\mathbf{D}_\Delta(\Phi) = -\Phi(\mathbf{Y}_n) \, , \quad (53)$$

whose solution $\Delta_n$ is used to update the next iterate as $\mathbf{Y}_{n+1} = \mathbf{Y}_n + \Delta_n$.

The perturbation matrix $\Delta$ appears linearly in eq. (52), as does $\Delta^\mathsf{T}$, and each is multiplied on both sides by different matrix coefficients. A system of equations of this form is known as a *general linear system.* It is also known that in order to solve it by standard LU decomposition methods it must be converted into the form

$$\mathbf{K}_n \, col(\Delta) = -\, col(\Phi(\mathbf{Y}_n)) \, . \tag{54}$$

Here, $col(\cdot)$ is an operator that stacks the columns of its argument matrix into a single column vector. The construction of the matrix $\mathbf{K}_n \in \mathbf{R}^{3N \times 3N}$ is based on the Kronecker product $\otimes$, and the relation $col(\mathbf{LMN}^\mathsf{T}) = (\mathbf{N} \otimes \mathbf{L}) \, col(\mathbf{M})$, which holds for arbitrary conformable matrices $\mathbf{L}$, $\mathbf{M}$ and $\mathbf{N}$ [11]. On applying the $col$ operator to eq. (52) and using this relation, we find that

$$\begin{aligned} col(\mathbf{D}_\Delta(\Phi)) = \, & [(\mathbf{Y}_n^\mathsf{T} \mathbf{Y}_n + \mathbf{X}^\mathsf{T}\mathbf{X}) \otimes \mathbf{I}_N] \, col(\Delta) + \\ & + [\mathbf{Y}_n^\mathsf{T} \otimes \mathbf{Y}_n] \, col(\Delta^\mathsf{T}) + ([\mathbf{I}_3 \otimes (\mathbf{Y}_n \mathbf{Y}_n^\mathsf{T} - \mathbf{C})] \, col(\Delta)) \, . \end{aligned} \tag{55}$$

To write this as a single linear system in $col(\Delta)$, let $\mathbf{E}_{ij} \in \mathbf{R}^{N \times 3}$ be an elementary matrix, which has a one in the $ij$-th position and zeros elsewhere, and let $\mathbf{P} \in \mathbf{R}^{3N \times 3N}$ be a permutation matrix whose columns are $\mathbf{P} = [col(\mathbf{E}_{11}), \ldots, col(\mathbf{E}_{13}), col(\mathbf{E}_{21}), \ldots, col(\mathbf{E}_{N3})]$. It follows that $col(\Delta^\mathsf{T}) = \mathbf{P} \, col(\Delta)$ (see e.g. [11]), and hence $\mathbf{K}_n$ is given by

$$\mathbf{K}_n = [\mathbf{Y}_n^\mathsf{T}\mathbf{Y}_n + \mathbf{X}^\mathsf{T}\mathbf{X}] \otimes \mathbf{I}_N + [\mathbf{Y}_n^\mathsf{T} \otimes \mathbf{Y}_n] \mathbf{P} + \mathbf{I}_3 \otimes [\mathbf{Y}_n \mathbf{Y}_n^\mathsf{T} - \mathbf{C}] \, . \tag{56}$$

Once $col(\Delta_n) = -\mathbf{K}_n^{-1} \, col(\Phi(\mathbf{Y}_n))$ has been obtained by standard numer-ical methods, it is easily put back into the desired increment matrix $\Delta_n$.

The convergence of this iteration would be greatly improved if there were an inexpensive method of generating good starting matrices $\mathbf{Y}_0$. If $\|\mathbf{G_0}\|_1 \ll \bar{\gamma}$ (see eq. (46)), then we can assume $\mathbf{B} = \mathbf{0}$ to a good approximation, in which case the $\gamma_k$ will be close to the three largest eigenvalues of $\mathbf{C}$ and we can use $\mathbf{Y}_0 = \hat{\mathbf{Y}} \Gamma^{1/2}$, where the columns of $\hat{\mathbf{Y}}$ are the corresponding eigenvectors. If $\|\mathbf{C}\|_1 \ll \lambda_{min}(\mathbf{G_0})$, then we can simply set $\mathbf{Y}_0 = \mathbf{V}\mathbf{G_0}^{-1}$. More general procedures for finding good starting matrices, and searching the space for the global optimum, will be presented later on. The Newton-Kantorovitch method has the advantage of being robust and easy to implement; it has the disadvantage of being relatively inefficient since it requires the repeated solution of a $3N \times 3N$ system of linear equations.

## 5.2. CAYLEY PARAMETRIZATION OF THE INERTIAL EQUATIONS

In this section we construct a method for solving eqs. (20)-(21) for the eigenvalues $\gamma_1, \gamma_2, \gamma_3$ and eigenvectors $\mathbf{r}_1, \mathbf{r}_2, \mathbf{r}_3$ of the inertial tensor. This method is based on Cayley representation of an orthogonal matrix [16], which is possibly the simplest way to parametrize a rotation. Consider a real vector $\omega$ and the associated $3 \times 3$ skew-symmetric matrix:

$$\Omega = \begin{bmatrix} 0 & -\omega_3 & \omega_2 \\ \omega_3 & 0 & -\omega_1 \\ -\omega_2 & \omega_1 & 0 \end{bmatrix} . \tag{57}$$

The *Cayley transform* is a rational matrix function which maps $\Omega$ into an orthogonal matrix $\mathbf{R}$:

$$\mathbf{R} = (\mathbf{I} - \Omega)/(\mathbf{I} + \Omega) . \tag{58}$$

If $\|\omega\| = \sqrt{\omega_1^2 + \omega_2^2 + \omega_3^2}$, the eigenvalues of $\mathbf{R}$ are given by $\{1, (1 + \iota\|\omega\|)/(1 - \iota\|\omega\|), (1 - \iota\|\omega\|)/(1 + \iota\|\omega\|)\}$, which shows that $\mathbf{R}$ is a proper orthogonal matrix, or rotation. The mapping $\omega \to \mathbf{R}$ is an isomorphism whose inverse is given by

$$\Omega = (\mathbf{I} - \mathbf{R})/(\mathbf{I} + \mathbf{R}) \tag{59}$$

whenever $\mathbf{R}$ is a rotation. By using the Cayley parametrization, we (a) eliminate all six constraints in eq. (21) and hence reduce the total number of variables from twelve to six, and (b) have an efficient way of com-puting the vector $\omega$ for any given rotation matrix (cf. eq. (59)).

In constructing $\mathbf{R}$ from the three eigenvectors $\mathbf{r}_1, \mathbf{r}_2, \mathbf{r}_3$, their arbitrary orientation causes an ambiguity in the values of $\omega$. There are eight possible orientations, but four of these are eliminated since they result in $det(\mathbf{R}) = -1$. Among the remaining four, we shall always choose the one that maximizes $tr(\mathbf{R})$, because (a) it results in the smallest magnitude of $\|\omega\|$, and (b) improves the condition number of the matrix division in eq. (59).

Appending the three rotational parameters $\omega_1, \omega_2, \omega_3$ to the three eigenvalues $\gamma_1, \gamma_2, \gamma_3$, we form a vector of parameters $\mathbf{p} \in \mathbf{R}^6$, and then write the system of eqs. (20) in vector form

$$\mathbf{g}(\mathbf{p}) = \begin{bmatrix} g_{11}(\mathbf{p}) & g_{22}(\mathbf{p}) & g_{33}(\mathbf{p}) & g_{12}(\mathbf{p}) & g_{13}(\mathbf{p}) & g_{23}(\mathbf{p}) \end{bmatrix}^\mathsf{T} = \mathbf{0} . \qquad (60)$$

Starting from a parameter vector $\mathbf{p}_0$, Newton's method for the system of equations (60) is given by $\mathbf{p}_{n+1} = \mathbf{p}_n - \mathbf{J}_n^{-1} \cdot \mathbf{g}_n$, where $\mathbf{g}_n = \mathbf{g}(\mathbf{p}_n)$ is the function vector and $\mathbf{J}_n = \partial \mathbf{g}/\partial \mathbf{p}|_{\mathbf{p}_n}$ is the Jacobian matrix evaluated at $\mathbf{p}_n$.

The partial derivatives of eq. (60) with respect to the Cayley parameters $\omega_1, \omega_2, \omega_3$ are computed as follows. Let $\Delta_m = \partial \Omega / \partial \omega_m$, $m = 1, 2, 3$, be the three directional derivative matrices of $\Omega$, i.e.

$$\Delta_1 = \begin{bmatrix} 0 & 0 & 0 \\ 0 & 0 & -1 \\ 0 & 1 & 0 \end{bmatrix} \quad \Delta_2 = \begin{bmatrix} 0 & 0 & 1 \\ 0 & 0 & 0 \\ -1 & 0 & 0 \end{bmatrix} \quad \Delta_3 = \begin{bmatrix} 0 & -1 & 0 \\ 1 & 0 & 0 \\ 0 & 0 & 0 \end{bmatrix}. \quad (61)$$

Upon differentiating the expression $(\mathbf{I} + \Omega)(\mathbf{I} + \Omega)^{-1} = \mathbf{I}$ with respect to $\omega_m$, we get $\Delta_m (\mathbf{I} + \Omega)^{-1} + (\mathbf{I} + \Omega)[\partial / \partial \omega_m (\mathbf{I} + \Omega)^{-1}] = \mathbf{0}$, or

$$\frac{\partial}{\partial \omega_m}(\mathbf{I} + \Omega)^{-1} = -(\mathbf{I} + \Omega)^{-1} \Delta_m (\mathbf{I} + \Omega)^{-1}. \quad (62)$$

We now substitute eq. (62) into $\partial \mathbf{R} / \partial \omega_m = (\mathbf{I} - \Omega)[\partial / \partial \omega_m (\mathbf{I} + \Omega)^{-1}] - \Delta_m (\mathbf{I} + \Omega)^{-1}$, simplify it and obtain

$$\frac{\partial \mathbf{R}}{\partial \omega_m} = -2(\mathbf{I} + \Omega)^{-1} \Delta_m (\mathbf{I} + \Omega)^{-1}. \quad (63)$$

The three directional derivative matrices, $\partial \mathbf{R} / \partial \omega_m$, contain all nine derivative vectors $\partial \mathbf{r}_k / \partial \omega_m$ ($1 \leq k, m \leq 3$) and, as columns of a rotation matrix, they are considered *independent* of the $\gamma_k$ variables.

Since the vector $\mathbf{r}_k$ lies on the intersection curve $\mathbf{r}(\gamma_k)$ (by ***Lemma***

**3.**), one might ask why do we ignore dependence of $\mathbf{r}_k$ on $\gamma_k$? The answer is that when we determine three unit vectors during the global search that approximately satisfy the nonlinear system of equations (60), they are only approximately orthogonal. Therefore, the inverse Cayley transform (eq. (59)) results in approximately skew-symmetric matrix $\tilde{\Omega}$. When we skew-symmetrize it, according to $\bar{\Omega} = (\tilde{\Omega} - \tilde{\Omega}^\mathsf{T})/2$, in order to generate the three correct rotation parameters from which to start Newton's method, the corresponding rotation matrix $\bar{\mathbf{R}}$ has columns that do not lie any longer on the original intersection curves (generated by the $\gamma_k$) but are only in their vicinities. The result is that the terms containing $d\mathbf{r}(\gamma_k)/d\gamma_k\big|_{\mathbf{r}_k}$ can not be accurately calculated and hence they are omitted from eqs. (64) below. This simplification allows the orthonormality of the three vectors to be maintained, at the price of a somewhat less than quadratic rate of convergence.

The entries of this simplified Jacobian matrix take on the following form:

$$\begin{aligned}\frac{\partial g_{kk}}{\partial \gamma_k} &= 1 + 2\mathbf{r}_k^\mathsf{T}[\mathbf{G_0} + \mathbf{W}^\mathsf{T}([\gamma_k\mathbf{I} - \Sigma]^{-3})\mathbf{W}]\mathbf{r}_k & 1 \leq k \leq 3 \\ \frac{\partial g_{kl}}{\partial \gamma_k} &= -\mathbf{r}_k^\mathsf{T}[\mathbf{G_0} + \mathbf{W}^\mathsf{T}([\gamma_k\mathbf{I} - \Sigma]^{-2}[\gamma_l\mathbf{I} - \Sigma]^{-1})\mathbf{W}]\mathbf{r}_l & 1 \leq k < l \leq 3\end{aligned} \quad (64)$$

Finally, by eqs. (18)-(19), $\mathbf{S}(\gamma_k, \gamma_l) \to \mathbf{S}(\gamma_k)$ as $\gamma_l \to \gamma_k$, so that

$$\frac{\partial g_{kl}}{\partial \omega_m} = \left[\frac{\partial \mathbf{r}_k}{\partial \omega_m}\right]^T \mathbf{S}(\gamma_k, \gamma_l) \mathbf{r}_l + \left[\frac{\partial \mathbf{r}_l}{\partial \omega_m}\right]^T \mathbf{S}(\gamma_k, \gamma_l) \mathbf{r}_k \qquad 1 \le k \le l \le 3 \ . \quad (65)$$

The rest of the entries are equal to zero.

## 5.3. ANGULAR PARAMETRIZATION OF THE INERTIAL EQUATIONS

An alternative approach to reducing eqs. (20)-(21) to a system of six equations in six unknowns is to use as the parameters the three eigenvalues $\gamma_k$ together with the three angles $\psi_k$ defined in eq. (30), one for each intersection curve $\mathbf{r}(\gamma_k)$. For any three values of $\gamma_k$ (satisfying the intersection condition (27)) and any three values of $\psi_k \in [0, 2\pi)$, the three unit vectors $\mathbf{r}(\gamma_k, \psi_k)$ satisfy the first three eqs. (20) and the three normalizing conditions in eqs. (21), but are not necessarily orthogonal. Conversely, for each column vector $\mathbf{r}_k$ on the intersection curve $\mathbf{r}(\gamma_k)$ we find the corresponding angular parameter $\psi_k$ by means of eqs. (29)-(31).

If we denote the parameter vector $\mathbf{q} = \begin{bmatrix} \gamma_1 & \psi_1 & \gamma_2 & \psi_2 & \gamma_3 & \psi_3 \end{bmatrix}^T$, then the remaining three eqs. (20) can be written as

$$h_1(\gamma_1, \psi_1, \gamma_2, \psi_2) = \mathbf{r}^T(\gamma_1, \psi_1)\, \mathbf{S}(\gamma_1, \gamma_2)\, \mathbf{r}(\gamma_2, \psi_2) = 0$$
$$h_3(\gamma_2, \psi_2, \gamma_3, \psi_3) = \mathbf{r}^T(\gamma_2, \psi_2)\, \mathbf{S}(\gamma_2, \gamma_3)\, \mathbf{r}(\gamma_3, \psi_3) = 0 \qquad (66)$$
$$h_5(\gamma_3, \psi_3, \gamma_1, \psi_1) = \mathbf{r}^T(\gamma_1, \psi_1)\, \mathbf{S}(\gamma_1, \gamma_3)\, \mathbf{r}(\gamma_3, \psi_3) = 0\ ,$$

while the three orthogonality conditions in eqs. (21) become

$$h_2(\gamma_1, \psi_1, \gamma_2, \psi_2) = \mathbf{r}^T(\gamma_1, \psi_1)\, \mathbf{r}(\gamma_2, \psi_2) = 0$$
$$h_4(\gamma_2, \psi_2, \gamma_3, \psi_3) = \mathbf{r}^T(\gamma_2, \psi_2)\, \mathbf{r}(\gamma_3, \psi_3) = 0 \qquad (67)$$
$$h_6(\gamma_3, \psi_3, \gamma_1, \psi_1) = \mathbf{r}^T(\gamma_1, \psi_1)\, \mathbf{r}(\gamma_3, \psi_3) = 0\ .$$

This completes the description of the system of equations.

We now illustrate the computation of the entries of the Jacobian matrix by presenting only the partial derivatives for the functions $h_1$ and $h_2$ with respect to $\gamma_1$ and $\psi_1$. The remaining derivatives have the same form apart from index permutation. The partial derivatives with respect to $\psi_1$ are straightforward:

$$\frac{\partial h_1}{\partial \psi_1} = \frac{\partial \mathbf{r}^T(\gamma_1, \psi_1)}{\partial \psi_1} \mathbf{S}(\gamma_1, \gamma_2) \, \mathbf{r}(\gamma_2, \psi_2)$$
$$\frac{\partial h_2}{\partial \psi_1} = \frac{\partial \mathbf{r}^T(\gamma_1, \psi_1)}{\partial \psi_1} \, \mathbf{r}(\gamma_2, \psi_2) \quad ,$$
(68)

because $\mathbf{r}(\gamma_2, \psi_2)$ does not depend on $\psi_1$. The partial derivative of the curve $\mathbf{r}(\gamma_1, \psi_1)$ with respect to $\psi_1$ is obtained via eqs. (29)-(31) as

$$\frac{\partial \mathbf{r}(\gamma_1, \psi_1)}{\partial \psi_1} = \mathbf{U}(\gamma_1) \begin{bmatrix} \partial \xi_1/\partial \psi_1 \\ \partial \xi_2/\partial \psi_1 \\ \partial \xi_3/\partial \psi_1 \end{bmatrix} = \mathbf{U}(\gamma_1) \begin{bmatrix} -\beta_1 \sin(\psi_1) \\ \beta_2 \cos(\psi_1) \\ \dfrac{(\beta_2 - \beta_1) \cos(\psi_1) \sin(\psi_1)}{\xi_3} \end{bmatrix} \quad (69)$$

if, for example, the axis of symmetry is $\mathbf{u}_3$.

The partial derivatives of $h_1$ and $h_2$ with respect to the eigenvalue $\gamma_1$ are more complicated. They depend on the derivatives of the spectral decomposition of the matrix $\mathbf{Q}(\gamma) = \mathbf{S}(\gamma)/\gamma$ (see section **4.1**). The derivatives of the eigenvalues and eigenvectors of a matrix-valued function of a scalar argument have been rediscovered many times, see e.g. [12], [17] and references therein. In physics, the eigenvalue derivatives are implicit in the Hellmann-Feynman theorem [18], but the eigenvector derivatives appear to be much less widely known. For this reason we include a brief derivation of these derivatives for the case of a symmetric (or hermitian) matrix.

Assuming that $\mathbf{Q}(\gamma) = \mathbf{U}(\gamma)\Lambda(\gamma)\mathbf{U}^T(\gamma)$ is a continuously differentiable symmetric matrix function in some domain, its derivative with respect to $\gamma$ is $\mathbf{Q}' = \mathbf{U}'\Lambda\mathbf{U}^T + \mathbf{U}\Lambda'\mathbf{U}^T + \mathbf{U}\Lambda\mathbf{U}'^T$. This implies

$$\Lambda' + [\Lambda, \mathbf{M}] = \mathbf{U}^T \mathbf{Q}' \mathbf{U}, \tag{70}$$

where the square brackets denote the matrix commutator, and the matrix $\mathbf{M} = -\mathbf{U}^T\mathbf{U}'$ is variously known as the Cartan matrix, logarithmic derivative, or multiplicative derivative of $\mathbf{U}^T$ [19]. Since the commutator of any matrix with a diagonal one has zeros on its diagonal, it follows that

$$\Lambda' = \mathit{diag}(\mathbf{U}^T \mathbf{Q}' \mathbf{U}). \tag{71}$$

The matrix $\mathbf{U}^T\mathbf{Q}'\mathbf{U}$ is symmetric (though not diagonal, because $[\mathbf{Q}, \mathbf{Q}'] \neq \mathbf{0}$), which implies that the commutator $[\Lambda, \mathbf{M}]$ is symmetric too, and consequently $\mathbf{M}$ must be skew-symmetric. Its entries $m_{ij}$ are obtained by simply writing out the off-diagonal entries in eq. (70) and then solving the resulting (scalar) linear equations, namely

$$m_{ij} = -\mathbf{u}_i^T \mathbf{u}_j' = \frac{\mathbf{u}_i^T \mathbf{Q}' \mathbf{u}_j}{\lambda_i - \lambda_j} \qquad (i \neq j) . \tag{72}$$

Finally, the eigenvector derivatives are computed from $\mathbf{U}' = -\mathbf{UM}$. When $\lambda_i = \lambda_j$ for a particular value of $\gamma$, a further analysis is required which we shall not pursue, because this case has not occurred in the applications of interest here.

Given the eigenvalue and eigenvector derivatives of $\mathbf{Q}(\gamma)$, the functions $h_1$ and $h_2$ are differentiated with respect to $\gamma_1$ as follows:

$$\frac{\partial h_1}{\partial \gamma_1} = \left( \frac{\partial \mathbf{r}^T(\gamma_1, \psi_1)}{\partial \gamma_1} \mathbf{S}(\gamma_1, \gamma_2) + \mathbf{r}^T(\gamma_1, \psi_1) \frac{\partial \mathbf{S}(\gamma_1, \gamma_2)}{\partial \gamma_1} \right) \mathbf{r}(\gamma_2, \psi_2)$$

$$\frac{\partial h_2}{\partial \gamma_1} = \frac{\partial \mathbf{r}^T(\gamma_1, \psi_1)}{\partial \gamma_1} \mathbf{r}(\gamma_2, \psi_2) . \tag{73}$$

The matrix derivative in eq. (73) is simply

$$\frac{\partial \mathbf{S}(\gamma_1, \gamma_2)}{\partial \gamma_1} = -\mathbf{W}^T ([\gamma_1 \mathbf{I}_N - \Sigma]^{-2} [\gamma_2 \mathbf{I}_N - \Sigma]^{-1}) \mathbf{W} , \tag{74}$$

while the partial derivative of the intersection curve with respect to $\gamma_1$

$$\frac{\partial \mathbf{r}(\gamma_1, \psi_1)}{\partial \gamma_1} = \frac{d\mathbf{U}(\gamma_1)}{d\gamma_1} \xi(\gamma_1, \psi_1) + \mathbf{U}(\gamma_1) \frac{\partial \xi(\gamma_1, \psi_1)}{\partial \gamma_1} \quad (75)$$

requires the matrix $d\mathbf{U}(\gamma_1)/d\gamma_1$ (given by eq. (72)), and the partial derivative vector $\partial \xi/\partial \gamma_1$. The last one depends on the ordinary derivatives of the projected semi-axes (see eq. (29)), for example,

$$\frac{d\beta_i}{d\gamma_1} = \begin{cases} \dfrac{(1-\lambda_2)\lambda_1' + (\lambda_1 - 1)\lambda_2'}{2\beta_i(\lambda_1 - \lambda_2)^2}, & \text{if } \lambda_2 < 1 \\ \dfrac{(1-\lambda_3)\lambda_1' + (\lambda_1 - 1)\lambda_3'}{2\beta_i(\lambda_1 - \lambda_3)^2}, & \text{if } \lambda_2 > 1 \end{cases} \quad (76)$$

which in turn depend on the eigenvalue derivatives $\lambda_m' = d\lambda_m/d\gamma_1$ (see eq. (71)). Thus, the expression for $\partial \xi/\partial \gamma_1$ is derived by using eqs. (29), (30) and (76).

# 6. The Search for the Global Minimum

## 6.1. AN OVERVIEW OF THE SEARCH STRATEGY

All three variations on Newton's method in section **5**. may diverge if they are not started sufficiently close to a solution. In this regard, the Newton-Kantorovitch method tends to be superior to the other two, but it requires the solution of a $3N \times 3N$ system of linear equations at each iteration, whereas the two based on the inertial equations require the solution of only a $6 \times 6$ system. Furthermore, due to the existence of several solutions, neither of these methods can be expected to find, on its own and with any degree of reliability, the *global* minimum of the STRAIN $f$, which is what interests us the most.

For this reason we have developed a search procedure that is capable of finding all solutions to inertial equations in the subspace defined by the bounds established in section **4**. The global minimum is then easily found by constructing the associated critical matrices via eq. (25), evaluating the function $f$ at each, and selecting the minimal one. Numerical evidence gathered over numerous problems of size $M + N < 200$ indicates that the number of such critical matrices rarely exceeds a dozen or so, and hence this search is generally quite efficient. The procedure as a whole consists of four stages.

In the first stage, the two $\gamma$-intervals defined by the bounds $\gamma_1 \in (max(\sigma_1, \lambda_{max}(\mathbf{G_0})), \gamma^+)$ and $\gamma_2, \gamma_3 \in (\gamma^-, \gamma_2^+)$ are discretized subject to a given curve length increment criterion with a deliberately increased discretization density in the vicinity of poles and bifurcation points. The outcome is two nonuniform discrete sequences, or meshes, one for $\gamma_1$, the other shared by both $\gamma_2, \gamma_3$. This stage of the algorithm takes a very small fraction of the overall time.

The main idea for the second stage comes from ***Lemma 3***, which shows that each discrete value of $\gamma_1$ in its range defines a unique intersection curve $\mathbf{r}(\gamma_1)$ which contains *all* solutions to the first inertial equation $\gamma_1 - \mathbf{r}_1^T \mathbf{S}(\gamma_1) \mathbf{r}_1 = 0$. By discretizing this curve, as described in section **6.2.**, we generate a discrete set of all inertial pairs $\{\gamma_1, \mathbf{r}_1\}$ at fixed $\gamma_1$. We shall do that for each curve defined by the $\gamma_1$-mesh and obtain a sample covering the portion of the unit sphere of interest. The discretization process for $\gamma$-intervals is described in section **6.3.** below.

The second stage systematically generates a discrete set of trial solutions over the $\gamma$-meshes, where each trial solution consists of three distinct inertial pairs. Each trial solution is tested to see if it satisfies the three bilinear inertial equations (20) to within a given threshold, and if so, it is appended to a list of initial solutions. These generation and elimination processes are described in section **6.4.** This stage of the algorithm may take a significant portion of the overall time, and is problem dependent.

The third and the fourth stage operate together as follows: each initial solution from the list is refined by means of one of the two versions of Newton's method derived in section **5**. In the case of convergence, the function $f$ is evaluated at this new critical matrix and inserted into an ordered list of the function values already found, provided it is different from them. These two final stages of the algorithm are proportional to the number of initial points generated in the second stage, and may also take significant portion of the overall time.

## 6.2. DISCRETIZATION OF THE INTERSECTION CURVES

For any given $\gamma_1$ in the mesh over the interval in eq. (35), let $\mathbf{r}(\gamma_1, \psi_1)$ be an intersection curve parametrized by $\psi_1$. We choose as our criterion for discretizing it a fixed spherical distance $\Delta\theta$ between successive points along the curve, namely

$$\mathbf{r}^\mathsf{T}(\gamma_1, \psi_1^k)\mathbf{r}(\gamma_1, \psi_1^{k+1}) = cos(\Delta\theta) \ . \tag{77}$$

Since we will be dealing in this section with a single fixed value of $\gamma_1$, we shall omit the subscript from $\psi_1^k$ and the dependence on $\gamma_1$ altogether, denoting the points on the intersection curve simply by $\mathbf{r}(\psi^k)$.

The problem then is to find a discrete sequence of parameters $\{\psi^k\}_{k=1}^m$ satisfying eq. (77).

Due to the double symmetry of the intersection curve (see section **4.1.**), it is sufficient to consider only the interval $[0, \pi/2]$, and find a subsequence of parameters $0 = \psi^0 < \psi^1 < \ldots < \psi^{m/4} = \pi/2$ satisfying (77). The complete sequence is then obtained by reflecting and copying this subsequence to obtain exactly $m$ points, as follows:

$$\begin{aligned} \{ \psi^0, \ldots, \psi^{m/4}, \quad & \pi - \psi^{m/4-1}, \ldots, \pi - \psi^0, \\ \pi + \psi^1, \ldots, \pi + \psi^{m/4}, \quad & 2\pi - \psi^{m/4-1}, \ldots, 2\pi - \psi^1 \} \end{aligned} \tag{78}$$

Suppose we have found points $\psi^0, \ldots, \psi^k$ satisfying eq. (77), and such that $\psi^k < \pi/2$. Finding the next point means solving the equation

$$\mathbf{r}^T(\psi^k)\mathbf{r}(\psi^k + \delta\psi) = cos(\Delta\theta) \tag{79}$$

for the increment $\delta\psi$. While it is certainly possible to solve eq. (79) at each step to full working precision, the final parameter will in general exceed $\pi/2$. If the last parameter is reset to $\pi/2$, then the last spherical distance may significantly differ from others. A slight redistribution of the points, on the other hand, would accommodate for this nonunifor-mity but then some of the work in getting the exact solution would have been wasted. A way out of this situation comes from the following argument. Consider the derivative of the arc length of the intersection curve with respect to the parameter $\psi$, i.e.

$$s'(\psi) = \sqrt{\frac{\beta_i^2 \sin^2(\psi) - \beta_j^2(\beta_i^2 - \cos^2(\psi))}{1 - \beta_i^2 \cos^2(\psi) - \beta_j^2 \sin^2(\psi)}}, \qquad (80)$$

which is obtained from $s'(\psi) = \sqrt{(\xi_1'(\psi))^2 + (\xi_2'(\psi))^2 + (\xi_3'(\psi))^2}$ and derivatives in eq. (30) and where $0 < \beta_i, \beta_j < 1$ are the semi-axes of the projected ellipse. A quick analysis shows that $s'(\psi)$ is a periodic function with period $\pi$, which is monotone on the intervals $[0, \pi/2]$ and $[\pi/2, \pi]$, but with different monotonicity on latter due to the reflection $s'(\psi) = s'(\pi - \psi)$ in $\pi/2$. Thus, $s'(\psi)$ attains its extremal values $min(\beta_i, \beta_j)$ and $max(\beta_i, \beta_j)$ at the points $0$ and $\pi/2$. Since the arc length between any two points on a continuous spherical curve is bounded from below by their spherical distance, we see that

$$\begin{aligned} \Delta\theta &= acos(\mathbf{r}^\mathsf{T}(\psi^k)\mathbf{r}(\psi^k + \delta\psi)) \\ &\leq \delta s^k = \int_{\psi^k}^{\psi^k + \delta\psi} s'(\psi)\, d\psi \leq max(\beta_i, \beta_j)\delta\psi \quad . \end{aligned} \qquad (81)$$

It follows that $\Delta\theta/max(\beta_i, \beta_j)$ is a lower bound on all the parametric increments $\delta\psi$ satisfying eq. (79). Locally, however, a better approximation is given by $\delta\psi \approx \Delta\theta/s'(\psi^k)$. The following discretization algorithm is based on this idea:

$$\begin{aligned}
set: \quad & \psi^0 = 0, \quad k = 0 \\
while \; & (\psi^k < \pi/2) \\
& \delta\psi^k = \Delta\theta/s'(\psi^k) \\
& \psi^{k+1} = \psi^k + \delta\psi^k \\
& k = k+1 \\
end &
\end{aligned} \quad (82)$$

Upon exiting the while loop, we set $m = 4k$ and then correct for any "overshoot" of $\pi/2$ by rescaling the sequence of increments $\delta\psi^k$ using the relative scale $\eta = \pi/(2\psi^{m/4}) < 1$. This means that a new sequence of angles is generated by $\psi^k = \psi^{k-1} + \delta\psi^{k-1}\eta$ for $k = 1, ..., m/4$.

This algorithm produces parameter sequences that approximately satisfy criterion (79). In addition, we make sure that it slightly oversamples the interval by starting the sequence from the point with the maximum derivative of the arc length, i.e. $\psi^0 = 0$ if $\beta_i \leq \beta_j$ or $\psi^0 = \pi/2$ if $\beta_i \geq \beta_j$, since then $s'(\psi^k) > s'(\psi^{k+1})$. In the last case, we first subtract $\pi/2$ from the computed subsequence before reflecting and copying as in eq. (78).

## 6.3. DISCRETIZATION OF THE γ INTERVALS

In order to discretize two γ intervals, $(max(\sigma_1, \lambda_{max}(\mathbf{G_0})), \gamma^+)$ and $(\gamma^-, \gamma_2^+)$, we need a suitable definition of distance between two consecutive nested curves $\mathbf{r}(\gamma^i)$ and $\mathbf{r}(\gamma^{i+1})$. We have considered a number of such definitions, e.g. the Hausdorff distance and several of its approximations, or the half-angle $\chi(\gamma)$ of the enveloping cone, but they resulted in either too many points, or are not sufficiently sensitive in the regions of fast variation of an intersection curve such as the tangency and bifurcation points.

It is these considerations that led us to consider the length of an intersection curve $0 \leq L(\gamma) \leq 2\pi$ as a quantity which is continuous, and inherits the behavior of both $\beta_i(\gamma)$ and $\beta_j(\gamma)$. Furthermore, if there is a bifurcation point $\gamma_b \in (max(\sigma_1, \lambda_{max}(\mathbf{G_0})), \gamma^+)$, then $L(\gamma)$ is a monotonically increasing function on the interval $\gamma \in (max(\sigma_1, \lambda_{max}(\mathbf{G_0})), \gamma_b)$, and a monotonically decreasing one on the interval $\gamma \in (\gamma_b, \gamma^+)$, with its maximum attained at the bifurcation point. In spite of the jump discontinuity in the values of $\beta_i(\gamma)$ and $\beta_j(\gamma)$ across $\gamma_b$, when we evaluate arc length derivatives in eq. (80) on both sides of $\gamma_b$, we find that $s'(\gamma_b^-, \psi) = s'(\gamma_b^+, \psi) = 1$ for all $\psi$, confirming that $max_\gamma(L(\gamma)) =$ $= L(\gamma_b) = 2\pi$. The derivative $L'(\gamma)$ is also well defined, except at the tangency points where $L(\gamma_t) = 0$ and $L'(\gamma_t) = \pm\infty$ due to a square root type singularity. At a bifurcation point we have

$$\lim_{\gamma \to \gamma_b^+} L'(\gamma) = \lim_{\gamma \to \gamma_b^-} L'(\gamma) = 0. \tag{83}$$

The length of an intersection curve is obtained by integrating the expression in Eq. (80) (which really is the partial derivative $s' = \partial s(\gamma, \psi)/\partial \psi$) and taking advantage of its double symmetry, i.e.

$$L(\gamma) = 4\int_0^{\pi/2} \frac{\partial s(\gamma, \psi)}{\partial \psi} d\psi. \tag{84}$$

For $dL/d\gamma$ we calculate the mixed partial derivative of $s(\gamma, \psi)$, so that

$$\frac{dL}{d\gamma} = 4\int_0^{\pi/2} \frac{\partial^2 s(\gamma, \psi)}{\partial \gamma \partial \psi} d\psi. \tag{85}$$

The interchange of integration and partial differentiation is permissible for any region not containing the points $\gamma_b$, $\gamma_t$, $\sigma_i$. Both integrals are evaluated using standard Romberg quadrature and the analytic expressions for $\partial s(\gamma, \psi)/\partial \psi$ (see eq. (80)) and $\partial^2 s(\gamma, \psi)/\partial \gamma \partial \psi$ (eqs. (80), (76) and (71)). The maximal number of interval halvings rarely exceeded 7 for the double precision accuracy.

In many ways the $\gamma$-discretization can be made similar to the $\psi$-discretization. This means that if we specify the length difference between two consecutive curves, namely $\Delta L = 2\pi/m$ for some $m > 1$, then instead of solving the nonlinear equation

$$L(\gamma^k + \delta\gamma) - L(\gamma^k) = \Delta L \tag{86}$$

exactly for the next increment $\delta\gamma$, we use an approximate version of eq. (86), namely

$$\gamma^{k+1} = \gamma^k + \frac{\Delta L}{L'(\gamma^k)} . \tag{87}$$

This will generate a decreasing sequence $\gamma^k$, providing it starts slightly to the left of $\gamma^+$ (to avoid division by $-\infty$), and ending slightly to the right of $max(\sigma_1, \gamma_b)$ (to avoid division by zero). In the case that $\sigma_1 < \gamma_b$ we continue generating the sequence of $\gamma^k$ using the same eq. (87), but with negative $\Delta L$, since the lengths are decreasing. As before, we start it slightly to the left of $\gamma_b$ and end it slightly to the right from $\sigma_1$. A small interval around a bifurcation point is discretized more densely in order to capture small or ill-conditioned intersection curves that are associated with solutions.

As mentioned earlier, when a bifurcation point is close to a pole, it often generates in its vicinity very thin intersection curves that appear as very narrow cuts on the unit sphere. Two kinds of problems occur in the presence of such an ill conditioned geometry: (a) the search process is not able to locate relevant initial solutions, and (b) Newton's method tends to diverge more frequently. We therefore increase the density of the $\gamma$-mesh in an inversely proportional manner to the eccentricity $min(\beta_i/\beta_j, \beta_j/\beta_i)$ in any region where this number drops below a given threshold. This mitigates the aforementioned difficulties.

## 6.4. ELIMINATION CRITERIA FOR THE GLOBAL SEARCH

We recall that for each $\gamma_1 \in (max(\sigma_1, \lambda_{max}(\mathbf{G_0})), \gamma^+]$ and for each $\psi_1 \in [0, 2\pi)$, the inertial pair $\gamma_1, \mathbf{r}_1 = \mathbf{r}(\gamma_1, \psi_1)$ satisfies the first inertial equation $\mathbf{r}_1^\mathsf{T}\mathbf{S}(\gamma_1)\mathbf{r}_1 = \gamma_1$. There are five more equations left to solve but only three variables, since we have fixed $\gamma_1$ and fixing the unit vector $\mathbf{r}_1$ eliminates *two* degrees of freedom. Two of these remaining variables are $\gamma_2$ and $\gamma_3$, while the third is chosen to be the angle $\phi$ in the plane orthogonal to $\mathbf{r}_1$ in which $\mathbf{r}_2$ and $\mathbf{r}_3$ must lie in order to form an orthonormal frame. Let $\mathbf{N} = [\mathbf{n}_1, \mathbf{n}_2]$ be any $3 \times 2$ matrix such that $[\mathbf{r}_1, \mathbf{n}_1, \mathbf{n}_2]$ forms a proper orthonormal frame. Then the vectors $\mathbf{r}_2, \mathbf{r}_3$ are expressed as

$$\begin{aligned} \mathbf{r}_2 &= \mathbf{n}_1 cos(\phi) + \mathbf{n}_2 sin(\phi) = \mathbf{N} \mathbf{t}_2 \\ \mathbf{r}_3 &= -\mathbf{n}_1 sin(\phi) + \mathbf{n}_2 cos(\phi) = \mathbf{N} \mathbf{t}_3 \,, \end{aligned} \quad (88)$$

where $\mathbf{t}_2^\top = [cos(\phi), sin(\phi)]$ and $\mathbf{t}_3^\top = [-sin(\phi), cos(\phi)]$ are evidently orthonormal plane vectors. We now rewrite the remaining two quadratic inertial equations in eqs. (20) in terms of this new angular variable $\phi$ as

$$\begin{aligned} \mathbf{t}_2^\top(\phi) \mathbf{N}^\top \mathbf{S}(\gamma_2) \mathbf{N} \mathbf{t}_2(\phi) &= \gamma_2 \\ \mathbf{t}_3^\top(\phi) \mathbf{N}^\top \mathbf{S}(\gamma_3) \mathbf{N} \mathbf{t}_3(\phi) &= \gamma_3 \end{aligned} \quad (89)$$

These two equations differ only by the subscripts on $\gamma$ and $\mathbf{t}$, and thus the solutions $\{\gamma, \phi\}$ to one of them are also the solutions to the other (recall that $\gamma_2, \gamma_3$ share the same domain as well as mesh). The same observation holds for the two bilinear equations $g_{12}, g_{13}$: they are identical, apart from the same two subscripts, namely

$$\begin{aligned} \mathbf{r}_1^\top \mathbf{S}(\gamma_1, \gamma_2) \mathbf{N} \mathbf{t}_2(\phi) &= 0 \\ \mathbf{r}_1^\top \mathbf{S}(\gamma_1, \gamma_3) \mathbf{N} \mathbf{t}_3(\phi) &= 0 \end{aligned} \quad (90)$$

These facts, which are a consequence of the formal symmetry of the sys-tem of eqs. (20), provide for considerable economy of effort. In the following analysis of these equations, therefore, we shall replace the subscript 2 or 3 with $t$ whenever the same argument applies equally to both of the equations in (89) or (90).

The $2 \times 2$ matrix $\mathbf{N}^T[\mathbf{S}(\gamma_t)/\gamma_t]\mathbf{N}$ in eq. (89) determines an ellipse, which is obtained by intersecting the ellipsoid (defined by $\mathbf{S}(\gamma_t)/\gamma_t$) with the central plane spanned by $[\mathbf{n}_1, \mathbf{n}_2]$. Therefore the solutions of eq. (89) are the points at which the intersection curves (between the ellipsoid and the unit sphere) intersect the plane $\mathbf{N}$. This is clearly a triple surface intersection problem, made considerably easier by the central symmetry of all three surfaces. Such an intersection exists if and only if the eigenvalues $\mu_{min}(\gamma_t)$ and $\mu_{max}(\gamma_t)$ of $\mathbf{N}^T[\mathbf{S}(\gamma_t)/\gamma_t]\mathbf{N}$ satisfy

$$\mu_{min}(\gamma_t) \leq 1 \leq \mu_{max}(\gamma_t) \tag{91}$$

Note also that, by the Rayleigh quotient argument, the inequalities $\lambda_3(\mathbf{S}(\gamma_t)/\gamma_t) \leq \mathbf{n}^T[\mathbf{S}(\gamma_t)/\gamma_t]\mathbf{n} \leq \lambda_1(\mathbf{S}(\gamma_t)/\gamma_t)$ hold for any unit vector $\mathbf{n}$, and in particular for $\mathbf{n} = \mathbf{Np}$ where $\mathbf{p}$ is any unit vector in the plane. This shows that $\lambda_3 \leq \mu_{min}$ and $\mu_{max} \leq \lambda_1$, so it is quite possible for the intersection curves not to intersect a given plane $\mathbf{N}$, and hence it is always necessary to check the inequalities in eq. (91) explicitly.

Let $\mathbf{P} = [\mathbf{p}_{min}, \mathbf{p}_{max}]$ be the eigenvectors of $\mathbf{N}^T[\mathbf{S}(\gamma_t)/\gamma_t]\mathbf{N}$ corresponding to the eigenvalues $\mu_{min}, \mu_{max}$. When the intersection condition (91) holds, the two exact intersection points for eq. (89) are given by

$$\mathbf{t}_\pm = \mathbf{P}\,\mathbf{m}_\pm, \qquad \mathbf{m}_\pm = \frac{1}{\sqrt{\mu_{max} - \mu_{min}}} \begin{bmatrix} \sqrt{1 - \mu_{min}} \\ \pm\sqrt{\mu_{max} - 1} \end{bmatrix}. \tag{92}$$

Their negatives, $-\mathbf{r}_\pm = -\mathbf{N}\,\mathbf{t}_\pm$, are also intersection points but, as we have seen, they are redundant in the sense that they give rise to the same critical matrices of $f$, and hence need not be considered.

The above arguments can be restated as follows: for any fixed inertial pair $\{\gamma_1, \mathbf{r}_1\}$ and any $\gamma_t$ taken from its mesh, either there are no solutions (criterion (91) is not satisfied) or else two solutions $\mathbf{r}_\pm$ are calculated from eqs. (88) and (92). We now subject both solutions to the first test, based on the bilinear eq. (90) and an orthogonality threshold $\varepsilon_\perp$, i.e.

$$\left| \mathbf{r}_1^\mathsf{T} \mathbf{S}(\gamma_1, \gamma_t) \mathbf{r}_\pm \right| < \varepsilon_\perp \left\| \mathbf{S}(\gamma_1, \gamma_t) \right\|_1 . \tag{93}$$

This means that if an inertial pair $\{\gamma_t, \mathbf{r}_t\}$ satisfies this inequality, then it is an approximate solution to $g_{12}(\gamma_1, \mathbf{r}_1, \gamma_t, \mathbf{r}_t) = 0$ as well as an exact solution to the two quadratic eqs. $g_{11}(\gamma_1, \mathbf{r}_1) = 0$ and $g_{tt}(\gamma_t, \mathbf{r}_t) = 0$. Each *inertial pair* $\{\gamma_t, \mathbf{r}_t\}$ that satisfies this test is appended as $\{\gamma^l, \mathbf{r}^l\}$ to a temporary list $L_T$ (associated with a single $\{\gamma_1, \mathbf{r}_1\}$), which will be used in the next step of the search.

In this next step, two additional tests are applied to each *pair* of elements $\{\gamma^{l_1}, \mathbf{r}^{l_1}\}$, $\{\gamma^{l_2}, \mathbf{r}^{l_2}\}$ from the list $L_T$. First, we check that the pair of vectors is approximately orthogonal, i.e.

$$(\mathbf{r}^{l_1})^\mathsf{T} \mathbf{r}^{l_2} < \varepsilon_\perp . \tag{94}$$

Second, we check that each pair that passes test (94) also approximately satisfies the last inertial equation, i.e.

$$\left| \mathbf{r}_2^\mathsf{T} \mathbf{S}(\gamma_2, \gamma_3) \mathbf{r}_3 \right| < \varepsilon_\perp \left\| \mathbf{S}(\gamma_2, \gamma_3) \right\|_1 , \tag{95}$$

where $\gamma_2 = max(\gamma^{l_1}, \gamma^{l_2})$, $\gamma_3 = min(\gamma^{l_1}, \gamma^{l_2})$ and $\mathbf{r}_2$, $\mathbf{r}_3$ are the corresponding vectors from the list $L_T$. Each pair of inertial pairs that pass the last two tests is then appended, together with $\{\gamma_1, \mathbf{r}_1\}$, to a list $L_I = \{\gamma_1^i, \mathbf{r}_1^i; \gamma_2^i, \mathbf{r}_2^i; \gamma_3^i, \mathbf{r}_3^i\}$ of initial (approximate) solutions for Newton refinement.

In the third stage of the search, one of the two versions of Newton's method (presented in the sections **5.2** and **5.3**) for solving the inertial equations is run on each of the initial solutions in the list $L_I$. The performance obtained with both of these methods is comparable, but the method based on the Cayley parameters is a bit faster and easier to implement, and is the method used in the examples presented in the next section.

In our implementation Newton's method we have imposed a damping factor on the step size during the first few iterations to reduce the frequency of diverging runs. Each time this iteration converges to a solution, we convert the solution into the corresponding critical matrix $\mathbf{Y}$ via eq. (25), and then evaluate $f(\mathbf{Y})$ via eq. (15). If this value is not already on the final list $L_C$, it is inserted there. The critical matrix with the smallest value of the *STRAIN* is taken as the global minimum.

It may happen that the initial $\gamma$-discretizations and the choice of $\varepsilon_\perp$ do not result in any initial solution, or results in too many. Furthermore, it may happen that Newton's method does not converge for any given initial point, or converges too often to the same solution. These cases are controlled by starting with relatively large $\varepsilon_\perp$ and relatively course discretization density of $\gamma$-mesh. Then, either $\varepsilon_\perp$ is reduced, or $\gamma$-mesh refined, or both, to ensure finding global minimum without excessive amount of redundant computation.

## 7. Numerical Experiments

The code implementing the global search procedure described in this paper was developed using the interactive MATLAB® numerical linear algebra system, and then converted into an independent 'C' program for a significantly faster execution. We have evaluated the resulting program on two kinds of test problems, one using randomly generated coordinates, and the other generated from the coordinates of structures in the Protein Data Bank (PDB). The procedures used to generate these test problems, and the results obtained with them, are the subject of this section.

## 7.1. RANDOM TEST PROBLEMS

The procedure used to generate random test problems is as follows: First, the coordinates of $\mathbf{X} \in \mathbf{R}^{M \times 3}$ and $\mathbf{Y} \in \mathbf{R}^{N \times 3}$ are chosen with a uniform distribution from within two balls centered at $\mathbf{c}_X, \mathbf{c}_Y \in \mathbf{R}^3$ with radii $r_X, r_Y > 0$. Second, both sets of points are translated to the centroid of $\mathbf{X}$ by subtracting the vector $\bar{\mathbf{x}} = -M^{-1} \sum_i \mathbf{x}_i$. Third, the inertial tensor of the $\mathbf{X}$ coordinates is diagonalized as $\mathbf{X}^T\mathbf{X} = \mathbf{LZL}^T$ ($\mathbf{Z} = \textit{diag}(\zeta_1, \zeta_2, \zeta_3)$ and $\zeta_1 \geq \zeta_2 \geq \zeta_3$), and both the $\mathbf{X}$ and $\mathbf{Y}$ coordinates are rotated by $\mathbf{L}$ so as to obtain principal axis coordinates for $\mathbf{X}$. Finally, the coordinates of both $\mathbf{X}$ and $\mathbf{Y}$ are scaled by $\zeta_1^{-1/2}$, so as to equalize the norm of $\mathbf{X}^T\mathbf{X}$ across different problems.

The matrices of squared distances, i.e.

$$\begin{aligned}
\mathbf{D}^A &= [D_{ij}^A]_{i,j=1}^{M,M} = [\|\mathbf{x}_i - \mathbf{x}_j\|^2]_{i,j=1}^{M,M} \\
\mathbf{D}^B &= [D_{ij}^B]_{i,j=1}^{M,N} = [\|\mathbf{x}_i - \mathbf{y}_j\|^2]_{i,j=1}^{M,N}, \\
\mathbf{D}^C &= [D_{ij}^C]_{i,j=1}^{N,N} = [\|\mathbf{y}_i - \mathbf{y}_j\|^2]_{i,j=1}^{N,N}
\end{aligned} \qquad (96)$$

are then calculated from these coordinates, and their elements perturbed by uniformly distributed random numbers $1 + \vartheta_{ij}$ as follows:

$$\begin{aligned}
\tilde{\mathbf{D}}^B &= [\tilde{D}^B_{ij}]^{M,N}_{i,j=1} \equiv [(1+\vartheta_{ij})^2 D^B_{ij}]^{M,N}_{i,j=1} \quad (\vartheta_{ij} \in (-\varepsilon_B, \varepsilon_B)) \\
\tilde{\mathbf{D}}^C &= [\tilde{D}^C_{ij}]^{N,N}_{i,j=1} \equiv [(1+\vartheta_{ij})^2 D^C_{ij}]^{N,N}_{i,j=1} \quad (\vartheta_{ij} \in (-\varepsilon_C, \varepsilon_C))
\end{aligned} \quad (97)$$

We stress that these perturbed distance matrices can no longer be embedded in a three-dimensional Euclidean space.

In order to have the blocks of the corresponding Gram matrices $\mathbf{A}$, $\mathbf{B}$ and $\mathbf{C}$ refer to a common origin and at the same time ensure that $\mathbf{A} = \mathbf{X}\mathbf{X}^\top$ is fit exactly, it is necessary place the origin in the $\mathbf{X}$ set, preferably at its centroid as above. This can be done directly from the distances by assigning unit mass to all the $\mathbf{X}$ coordinates, and zero to the $\mathbf{Y}$ coordinates. Equation (4) then yields

$$\begin{aligned}
D^A_{0i} &= \frac{1}{M}\sum_{j=1}^{M} D^A_{ij} - \frac{1}{M^2}\sum_{1=j<k}^{M,M} D^A_{jk} \quad (i=1,\ldots,M) \\
\tilde{D}^C_{0i} &= \frac{1}{M}\sum_{j=1}^{M} \tilde{D}^B_{ij} - \frac{1}{M^2}\sum_{1=j<k}^{M,M} D^A_{jk} \quad (i=1,\ldots,N)
\end{aligned} \quad (98)$$

Finally, we set $\mathbf{A} = \mathbf{X}\mathbf{X}^\top$ and

$$\begin{aligned}
\mathbf{B} &= [(D^A_{0i} + \tilde{D}^C_{0j} - \tilde{D}^B_{ij})/2]^{M,N}_{i,j=1} \\
\mathbf{C} &= [(\tilde{D}^C_{0i} + \tilde{D}^C_{0j} - \tilde{D}^C_{ij})/2]^{N,N}_{i,j=1}
\end{aligned} \quad (99)$$

Thus the free parameters that define our test problems include the dimensions $M$ and $N$ of $\mathbf{X}$ and $\mathbf{Y}$, the ratio of distance between the centers of the spheres to the sum of their radii, $s_{XY} = \|\mathbf{c}_X - \mathbf{c}_Y\|/(r_X + r_Y)$, and the perturbation magnitudes $\varepsilon_B$ and $\varepsilon_C$. In the problems reported below, we set $\{M = 10, N = 50\}$, $\{M = N = 30\}$ and $\{M = 50, N = 10\}$ with $r_X = M^{1/3}$ and $r_Y = N^{1/3}$, so as to keep the density approximately constant in each case. In addition, we set $s_{XY} = 0$, $s_{XY} = 1/2$ and $s_{XY} = 1$ for each of these three cases, as well as $\varepsilon_B = \varepsilon_C = \varepsilon$ for three perturbation levels, $\varepsilon = 0.005, 0.05, 0.5$. We generate and solve ten random problems in each class, so that the total number of random problems is $10 \times 3 \times 3 \times 3 = 270$. All problems were initialized with an arc length increment of $\Delta L = \pi/6$ (for the sequence of intersection curves) and an orthogonality threshold of $\varepsilon_\perp = \sqrt{1/2}\, asin(\pi/6)$.

For each of the 27 classes of problems, we report the average over the 10 problems solved of: (a) the number of discretized points generated for $\gamma_1$ and for $\gamma_2$; (b) the total number of unit vectors $\mathbf{r}_1$ generated on the intersection curves for entire $\gamma_1$-mesh; (c) the number of initial points to which Newton's method was applied; (d) the number of points for which Newton's method converged; (e) the number of distinct solutions to the inertial equations located; (f) the minimum function value found; (g) the CPU time required for the entire procedure (in seconds). These statistics are summarized in Table 1.

**Table 1. Statistics for Random Test Problems (see text).**

| ε (%) | M, N | $s_{XY}$ | (a) | (b) | (c) | (d) | (e) | (f) | (g) |
|---|---|---|---|---|---|---|---|---|---|
| 0.5 | 50, 10 | 0.0 | 8, 29 | 58 | 5836 | 4959 | 1.0 | 0.0000 | 672 |
| 0.5 | 50, 10 | 0.5 | 22, 45 | 219 | 2234 | 1270 | 1.0 | 0.0002 | 315 |
| 0.5 | 50, 10 | 1.0 | 22, 68 | 199 | 4136 | 2010 | 1.0 | 0.0007 | 554 |
| 0.5 | 30, 30 | 0.0 | 25, 43 | 218 | 2563 | 623 | 1.0 | 0.0003 | 485 |
| 0.5 | 30, 30 | 0.5 | 21, 72 | 161 | 4026 | 931 | 1.1 | 0.0012 | 615 |
| 0.5 | 30, 30 | 1.0 | 21, 71 | 140 | 5001 | 641 | 1.1 | 0.0077 | 515 |
| 0.5 | 10, 50 | 0.0 | 22, 91 | 70 | 6532 | 743 | 3.2 | 0.0248 | 838 |
| 0.5 | 10, 50 | 0.5 | 21, 90 | 43 | 6755 | 725 | 4.2 | 0.0514 | 813 |
| 0.5 | 10, 50 | 1.0 | 21, 84 | 35 | 6817 | 874 | 2.4 | 0.1443 | 1018 |
| 5.0 | 50, 10 | 0.0 | 12, 31 | 65 | 5094 | 4144 | 1.0 | 0.0022 | 533 |
| 5.0 | 50, 10 | 0.5 | 25, 46 | 290 | 766 | 439 | 1.0 | 0.0207 | 127 |
| 5.0 | 50, 10 | 1.0 | 24, 69 | 253 | 2431 | 1180 | 1.1 | 0.0726 | 303 |
| 5.0 | 30, 30 | 0.0 | 33, 43 | 320 | 1984 | 514 | 1.0 | 0.0273 | 410 |
| 5.0 | 30, 30 | 0.5 | 24, 83 | 183 | 5997 | 650 | 1.2 | 0.1340 | 803 |
| 5.0 | 30, 30 | 1.0 | 24, 86 | 152 | 5587 | 881 | 1.2 | 0.8100 | 599 |
| 5.0 | 10, 50 | 0.0 | 24, 96 | 146 | 5828 | 430 | 3.2 | 2.1188 | 698 |
| 5.0 | 10, 50 | 0.5 | 24, 89 | 76 | 6123 | 880 | 2.9 | 7.7206 | 699 |
| 5.0 | 10, 50 | 1.0 | 24, 95 | 47 | 6246 | 731 | 2.8 | 10.038 | 821 |
| 50.0 | 50, 10 | 0.0 | 9, 31 | 67 | 6222 | 4630 | 1.0 | 0.2135 | 713 |
| 50.0 | 50, 10 | 0.5 | 22, 44 | 221 | 2440 | 1460 | 1.0 | 1.8616 | 325 |
| 50.0 | 50, 10 | 1.0 | 22, 66 | 219 | 2989 | 1333 | 1.0 | 7.0796 | 394 |
| 50.0 | 30, 30 | 0.0 | 35, 49 | 285 | 3235 | 728 | 1.0 | 5.2963 | 655 |
| 50.0 | 30, 30 | 0.5 | 21, 70 | 192 | 5177 | 537 | 1.1 | 18.166 | 696 |
| 50.0 | 30, 30 | 1.0 | 21, 69 | 136 | 5771 | 1542 | 1.2 | 82.830 | 780 |
| 50.0 | 10, 50 | 0.0 | 21, 89 | 78 | 6395 | 411 | 3.7 | 203.86 | 759 |
| 50.0 | 10, 50 | 0.5 | 21, 98 | 38 | 6892 | 1003 | 3.0 | 471.46 | 842 |
| 50.0 | 10, 50 | 1.0 | 21, 103 | 45 | 6850 | 782 | 3.5 | 1292.1 | 901 |

It may be seen that the number of critical matrices found (e) was usually just one, the exception being those problems in which **Y** contained more points (50) than **X** (10), where it averaged about four. The relatively small number of critical matrices satisfying bounds is of course what makes the global search feasible. We observed that the values of $\gamma_1$ and $\gamma_3$ at the solutions were generally quite close to $\gamma^+$ and $\gamma^-$, respectively. The values of $\gamma_2$, on the other hand, were sometimes significantly smaller than $\gamma_2^+$, particularly when a bifurcation point was present (as in most of the $s_{XY} = 1$ cases). We also computed the standard deviations for all the averages given in Table **I** (not shown). The standard deviations in the number of approximate solutions found (c) and the number of times Newton's method converged (d) were always comparable to the averages for these quantities, and since these calculations consumed most of the CPU time, the average times reported in (g) fluctuated widely even within a single class of problem.

Numerical experience so far indicates that, with a moderately small discretization step, the search procedure always succeeds in finding the global minimum of the *STRAIN*. The difficulties encountered in an earlier version of the code, which occurred when the global minimum lay in the immediate vicinity of a pole, tangency or bifurcation point, were eliminated by: a) improvements in the discretization criterion; b) increasing the mesh density in the vicinity of such points; c) decreasing the orthogonality threshold $\varepsilon_\perp$ and the discretization criteria for the $\gamma_1$ and $\gamma_2$ meshes whenever a run failed to find any critical points.

## 7.2. CHEMICAL TEST PROBLEMS

For our chemical test problems, we used two biologically representative systems. The first was the crystal structure of a complex of the enzyme dihydrofolate reductase (DHFR; entry 3DRC in the Brookhaven Protein Data Bank) with its inhibitor methotrexate (MTX) [20]. Only the nonhydrogen atoms of one of the two molecules in the asymmetric unit were used, which gave 33 MTX atoms for the unknown **Y** coordinates, 1262 DHFR atoms for the fixed **X** coordinates. The second was the crystal structure of a small protein, bovine pancreatic trypsin inhibitor (BPTI; entry 4PTI in the Protein Data Bank), once again without hydrogen atoms [21]. In this case the coordinates of the 284 backbone nitrogen, alpha-carbon, beta-carbon, carbonyl-carbon and carbonyl-oxygen atoms were used for the **X** set, while the remaining 170 sidechain atoms were used for the **Y** set.

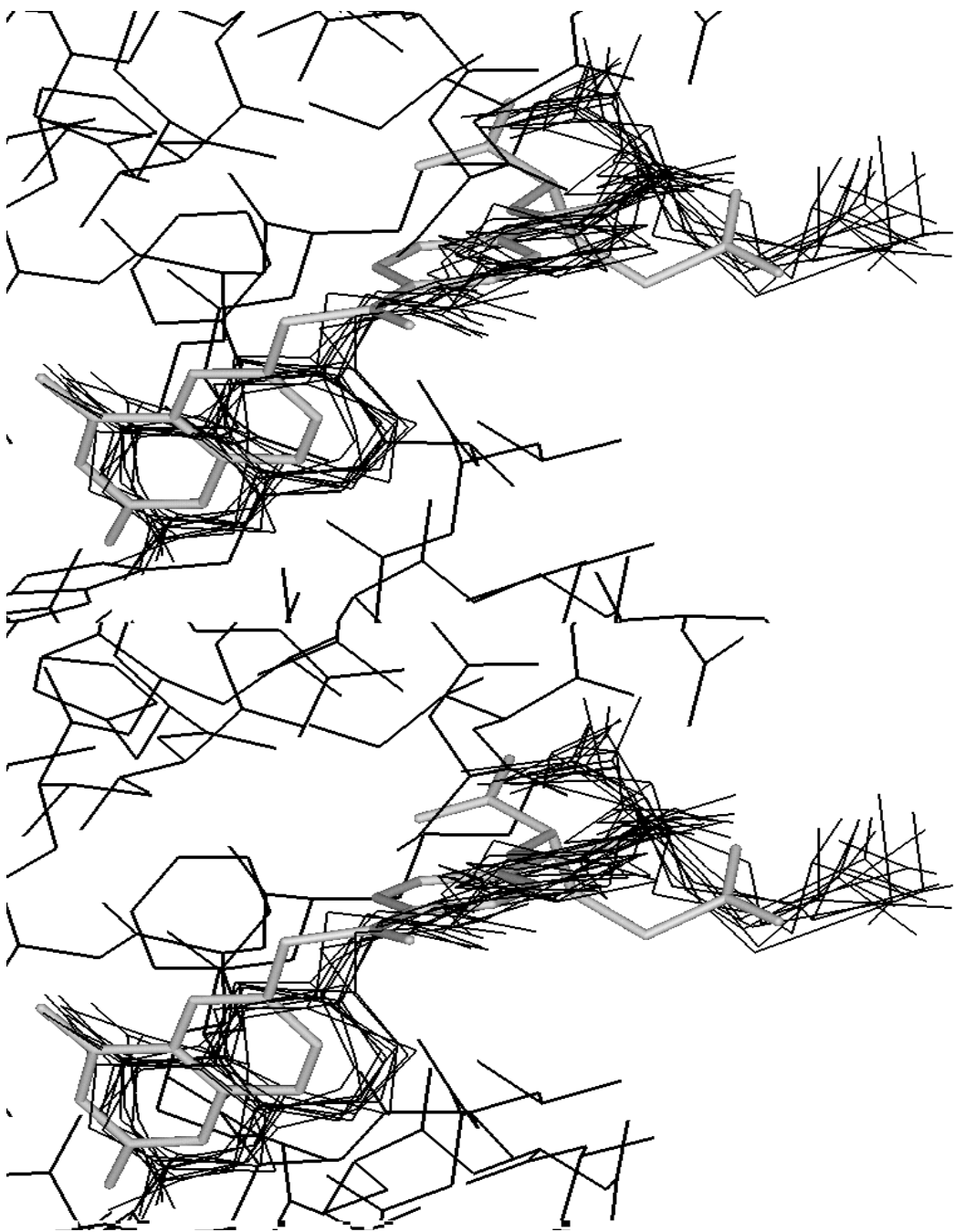

**Figure 5.** The results of DHFR/MTX test problems. The 10 MTX structures generated are drawn with medium lines, while the DHFR protein is drawn with thin lines, and the MTX crystal structure is a shaded stick drawing (to be viewed from the left with a stereoviewer).

In both cases the distance matrices were calculated from the coordinates in the crystal structures, and those distances involving the chosen **Y** coordinates were perturbed. The Gram matrices **B** and **C** were generated from these perturbed distances exactly as described for the random test problems above, but a somewhat different procedure was used to perturb the distances. First, if the distance $d$ was less than three Ångstroms, it was used unchanged. Second, if the distance $d$ was greater than three Ångstroms, it was replaced by a uniform random number between $max(3.0, d/\sqrt{2})$ and $d\sqrt{2}$. In the case of BPTI, this procedure generates random distance matrices like those that would be used as input to the *EMBED* algorithm in a homology modeling problem wherein only the backbone conformation of the protein was known. In the case of DHFR/MTX, the procedure generates random distance matrices like those that would be available if the structure of the protein was already known, and a number of ligand-protein and ligand-ligand interatomic contacts had been independently identified from NMR data.

The total number of random perturbations solved was 10 for both the DHFR/MTX and the BPTI backbone/sidechains problems. The DHFR/MTX trials ran smoothly, and although the number of approximate solu-tions found exceeded 6000 in all 10 trials, Newton's method converged better than 95% of the time to the same critical matrix every time, which is therefore almost certainly the unique global minimum. A single unique global minimum was also found in every case for BPTI, although Newton's method failed to converge much more often. This is due to the fact that the $\gamma_2$ was always close to the bifurcation point $\gamma_b$, and the angle between the curves at $\gamma_b$ was always very small (cf. section **4.1**). The time required for the BPTI problems ranged from 3 to 15 minutes.

Figures 5 and 6 show the results of the DHFR/MTX and BPTI test problems. Although a certain scatter is apparent in these coordinates, they are all very close to the crystal structure coordinates from which the unperturbed distances were obtained. This ability to correct for large (ca. 40%) but random perturbations in the distances is what makes the minimization criterion (9) and the algorithm a powerful tool for determination of structures. The figures also show that the perturbation procedure introduced some bias away from the crystal structure, which is expected when a simple uniform distribution is used for the perturbed distances [7]. How the procedure performs with the estimated distances that are available in actual applications will be the subject of a further study.

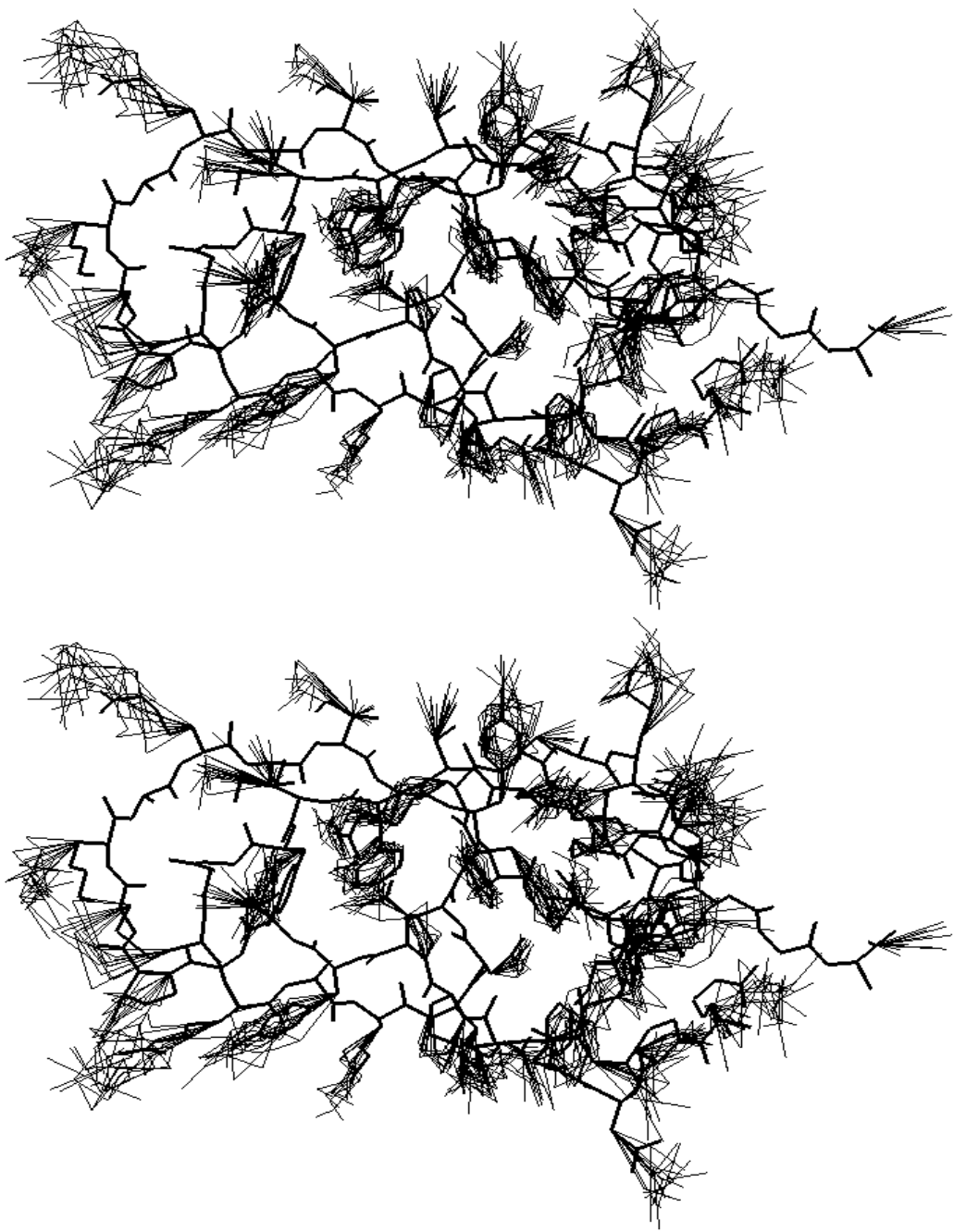

**Figure 6. The results of BPTI test problems treating the backbone atoms as fixed and the sidechain atoms as variable. The backbones of the computed structures have been superimposed on that of the crystal structure, which is drawn with a heavy line for comparison (see Fig. 5).**

# 8. Closing Comments

In this paper we have given a precise mathematical statement of the problem of embedding with a rigid substructure, derived necessary conditions for the global minimum, developed a search strategy for finding it, and evaluated the resulting algorithm on a realistic set of test problems. The rigid embedding problem, however, is only one step in a longer sequence of calculations by which one generates atomic coordinates that satisfy distance constraints. This means that in actual applications it will be necessary to have methods of incorporating fixed atom constraints into the other steps of the overall calculation. These include bound smoothing, the generation of approximate distances from which to build the metric matrix for embedding, and the optimization of the resulting coordinates.

It turns out that these steps are all relatively straightforward. Fixed atoms are built into the bound smoothing and distance generation steps simply by setting the distances among the fixed atoms to their known values. Similarly, the fixed atom constraints are incorporated into the optimization simply by eliminating their coordinates from the list of variables to be used. A suite of programs derived from the *DG-II* distance geometry program package [7] which incorporate these modifications will be described elsewhere, along with their application to chemically and biologically important problems [22].

## ACKNOWLEDGEMENTS

We thank Sven Hyberts for help with chemical test problems. This work was supported by NSF grants BIR-9511892 and MCB-9527181, by NIH grants GM-38221 and GM-47467, and by a grant from Molecular Simulations, Inc.